\documentstyle{article}

\renewcommand{\theequation}{\thesection.\arabic{equation}}
\newcommand{\sectie}[1]{\setcounter{equation}{0}\section{#1}}

\newcommand{\inlabel}[1]{\refstepcounter{equation} \label{#1}
\hspace{1.5ex} (\theequation)}

\newcommand{\N}{{\rm I\kern-.5ex N}}
\newcommand{\Z}{{\sf \vrule height 1.55ex depth-1.2ex
width.03em\kern-.11em Z
\kern-.9ex Z\kern-.11em\vrule height 0.3ex depth0ex width.03em}}
\newcommand{\Q}{{\rm\kern.2ex\vrule height1.55ex depth-.05ex
width.03em\kern-.7ex Q}}
\newcommand{\R}{{\rm I\kern-.5ex R}}
\newcommand{\Rvar}{{\rm I\kern-.5ex R}}
\newcommand{\C}{{\rm\kern.3ex\vrule height1.55ex depth-.05ex
width.03em\kern-.7ex C}}
\newcommand{\Cvar}{{\, \rm\kern.3ex\vrule height1.1ex depth-.05ex
width.03em\kern-.7ex C}}

\newcommand{\nabp}{\nab \hspace{-1.05ex}
\rule[.5ex]{.2ex}{.8ex}   \hspace{1.05ex}}
\newcommand{\sip}{\si \hspace{-0.92ex}
\rule[.15ex]{.1ex}{.55ex}   \hspace{0.92ex}}

\newcommand{\restr}{\hspace{0.15ex}\rule[-0.4ex]{.1ex}{1.5ex}\hspace{0.1ex}}

\newcommand{\spat}{\hspace{4ex}}

\newcommand{\flip}{ \chi }

\newcommand{\ahh}{\hat{A}
\hspace{-.55ex}\hat{\rule{0ex}{2.0ex}}\hspace{.55ex}}
\newcommand{\dhh}{\hat{\de}
\hspace{-.95ex}\hat{\rule{0ex}{2.05ex}}\hspace{.95ex}}

\newcommand{\ah}{\hat{A}}
\newcommand{\lambdah}{\hat{\lambda}}
\newcommand{\deh}{\hat{\Delta}}

\newcommand{\tauh}{\hat{\tau}}

\newcommand{\pih}{\hat{\pi}}

\newcommand{\nab}{\nabla}

\newcommand{\cI}{{\cal I}}
\newcommand{\cL}{{\cal L}}

\newcommand{\cN}{{\cal N}}
\newcommand{\cM}{{\cal M}}

\newcommand{\dual}{L^1(A)}
\newcommand{\duals}{L^1_*(A)}

\newcommand{\ot}{\otimes}

\newcommand{\la}{\Lambda}

\newcommand{\om}{\omega}
\newcommand{\io}{\iota}
\newcommand{\vfi}{\varphi}
\newcommand{\vep}{\varepsilon}

\newcommand{\al}{\alpha}
\newcommand{\be}{\beta}
\newcommand{\ga}{\Gamma}
\newcommand{\sde}{\delta}
\newcommand{\de}{\Delta}

\newcommand{\th}{\theta}

\newcommand{\si}{\sigma}
\newcommand{\Mfi}{{\cal M}_{\vfi}}
\newcommand{\Nfi}{{\cal N}_{\vfi}}

\newcommand{\Npsi}{{\cal N}_{\psi}}

\newcommand{\ovt}{\, \bar{\otimes}\, }
\newcommand{\siup}{\sip^{\text{\tiny u}}}
\newcommand{\siu}{\si^{\text{\tiny u}}}
\newcommand{\psiu}{\psi_{\text{\tiny u}}}
\newcommand{\vfiu}{\vfi_{\text{\tiny u}}}
\newcommand{\Ru}{R_{\text{\tiny u}}}
\newcommand{\tauu}{\tau^{\text{\tiny u}}}
\newcommand{\lambdau}{\lambda_{\text{\tiny u}}}
\newcommand{\lambdahu}{\hat{\lambda}_{\text{\tiny u}}}
\newcommand{\lau}{\la_{\text{\tiny u}}}
\newcommand{\gau}{\ga_{\text{\tiny u}}}
\newcommand{\deu}{\de_{\text{\tiny u}}}
\newcommand{\sdeu}{\sde_{\text{\tiny u}}}
\newcommand{\Su}{S_{\text{\tiny u}}}
\newcommand{\vepu}{\vep_{\text{\tiny u}}}
\newcommand{\taudu}{\tauh^{\text{\tiny u}}}
\newcommand{\ahu}{\hat{A}_{\text{\tiny u}}}
\newcommand{\dehu}{\deh_{\text{\tiny u}}}
\newcommand{\au}{A_{\text{\tiny u}}}
\newcommand{\alu}{\al_{\text{\tiny u}}}
\newcommand{\beu}{\be_{\text{\tiny u}}}
\newcommand{\alr}{\al_{\text{\tiny r}}}
\newcommand{\ber}{\be_{\text{\tiny r}}}
\newcommand{\gammau}{\gamma_{\text{\tiny u}}}

\newcommand{\cU}{{\cal U}}
\newcommand{\cV}{{\cal V}}

\newcommand{\cT}{{\cal T}}

\newcommand{\cA}{{\cal A}}

\newcommand{\text}[1]{\mbox{#1}}

\newcommand{\cst}{\text{C}$\hspace{0.1mm}^*$}
\newcommand{\wst}{\text{W}$\hspace{0.1mm}^*$}

\newcommand{\qed}{\ \hfill \rule{2mm}{2mm}}
\newenvironment{demo}{\medskip\noindent\bf Proof :\ \
\rm}{\qed\bigskip\par }

\newtheorem{definition}{Definition}[section]
\newtheorem{proposition}[definition]{Proposition}
\newtheorem{lemma}[definition]{Lemma}
\newtheorem{corollary}[definition]{Corollary}
\newtheorem{remark}[definition]{Remark}
\newtheorem{theorem}[definition]{Theorem}

\newtheorem{notation}[definition]{Notation}
\newtheorem{result}[definition]{Result}

\newtheorem{terminology}[definition]{Terminology}

\begin{document}

\begin{center}
\LARGE\bf Locally compact quantum groups in the universal setting \end{center}

\bigskip\bigskip

\begin{center}
\rm Johan Kustermans

Department of Mathematics

University College Cork

Western Road

Cork

Ireland

\medskip

e-mail : johank@ucc.ie

\bigskip\bigskip

\bf January 1999 \rm
\end{center}

\bigskip

\begin{abstract}
\noindent In this paper we associate to every reduced \cst-algebraic quantum group $(A,\de)$ (as defined in \cite{J-V}) a universal \cst-algebraic quantum group $(\au,\deu)$. We fine tune a proof of Kirchberg to show that every $^*$-representation of a modified $L^1$-space is generated by a unitary corepresentation. By taking the universal enveloping \cst-algebra of a dense sub $^*$-algebra of $A$ we arrive at the \cst-algebra $\au$. We show that this \cst-algebra $\au$ carries a quantum group structure which is as rich as its reduced companion. 
\end{abstract} 

\bigskip\medskip

\section*{Introduction}

In 1977, S.L. Woronowicz proposed the use of the \cst-language to axiomatize quantizations of locally compact quantum groups. This approach was very successful in the compact case (\cite{Wor7},\cite{Wor4},\cite{VD1}) and the discrete case (\cite{PW},\cite{VD3},\cite{ER}). In both cases the existence of the Haar weights could be proven from a simple set of axioms. The situation for the general non-compact however is less satisfactory. At present, there is still no general definition for a  locally compact quantum group in which the existence of the Haar weights is not one of the axioms of the proposed definition.

\medskip

The first attempt to axiomatize locally quantum groups aimed at enlarging the category of locally compact quantum groups in such a way that it contains locally compact groups and the reduced group \cst-algebras. A complete solution for this problem was found independently by M. Enock \& J.-M. Schwarz and by Kac \& Vainermann (see \cite{E} for a detailed account). The resulting objects are called Kac algebras and their definition was formulated in the von Neumann algebra framework. For quite a time, the main disadvantage of this theory lay in the fact that there was a lack of interesting examples aside from the groups and group duals.

\smallskip

S.L. Woronowicz constructed in \cite{Wor1} quantum $SU(2)$, an object which has all the right properties to be called a compact quantum group but does not fit into the framework of Kac algebras. In subsequent papers (\cite{Wor7},\cite{Wor4}), S.L. Woronowicz developed the axiom scheme for compact quantum groups. In contrast to the Kac algebra theory, quantum SU(2) fitted into this category of compact quantum groups.

The main difference between compact Kac algebras and compact quantum groups according to Woronowicz lies in the fact that the antipode of the Kac algebra is an automorphism  while in the approach of Woronowicz, it can be unbounded (as it is the case for quantum $SU(2)$).

\smallskip

It was Kirchberg (\cite{Kir}) who proposed a generalized axiom scheme for quantum groups in which the  antipode was unbounded but in which the antipode could be decomposed in an automorphism and an unbounded operator generated by a one-parameter group. This decomposition is called the polar decomposition of the antipode. This polar decomposition appeared for the compact case in \cite{Wor4}.  The general case was treated in the von Neumann algebra setting in \cite{Mas-Nak} by Masuda \& Nakagami. The main problem of their proposed definition of a quantum group lies in the complexity of the axioms.

In \cite{J-V}, the author and S.Vaes propose a relatively simple definition of a locally compact quantum group in its reduced form, i.e. in the form for which the Haar weights are faithful. We start of with a \cst-algebra with a comultiplication satisfying some density conditions and assume the existence of a faithful left invariant weight and a right invariant weight satisfying some kind of KMS condition. From these axioms, we are able to construct the antipode and its polar decomposition, prove the uniqueness of the Haar weights and construct the modular element. In short, we prove that the polar decomposition of the antipode is a consequence of some natural KMS assumptions on the Haar weights.

\smallskip

All the general axiom schemes considered above (except for the compact quantum groups by S.L. Woronowicz) are stated in the reduced setting. In this setting the von Neumann algebra approach and the \cst-algebra approach are equivalent and are in fact nothing else but two different ways a quantum group can present itself.

A quantum group can present itself in a third natural way, the universal way. In this case, one starts with a reduced locally compact quantum group $(A,\de)$. Then one considers a natural dense sub $^*$-algebra $B$ inside $A$ and proves that this $^*$-algebra has a universal enveloping \cst-algebra $\au$. The aim of this paper is to show that this universal \cst-algebra $\au$ carries a quantum group structure which is as rich on the analytical level as the reduced companion $A$. However, in this universal setting, the Haar weights do not have to be faithful. In return,we get the existence of a bounded counit.

\smallskip

In considering the universal dual of a quantum group (as opposed to its reduced dual), one is able to get a bijection between non-degenerate $^*$-representations of this universal dual and the unitary corepresentations of the original quantum group. This difference between reduced and universal duals is a mere generalization of the difference between the reduced and universal group \cst-algebras of a locally compact  group.

\medskip

The paper is organized as follows. In the first section, we fine tune a proof of Kirchberg to prove that every $^*$-representation of a modified $L^1$-space of $A$ is generated by a unitary corepresentation of $(A,\de)$. In the second section we introduce the comultiplication $\deu$ and counit $\vepu$ in the standard way (see \cite{PW}). We also construct the universal corepresentation of $(\au,\deu)$. The third section revolves around a procedure to lift automorphism on $A$ commuting with $\de$ from the reduced to the universal level. The Haar weights of $(\au,\deu)$ are introduced in section 4. In section 5, we construct the antipode and its polar decomposition. In section 6, we lift the modular element from the reduced to the universal level.

\bigskip\medskip

\section*{Notations and conventions}

For any subset $X$ of a Banach space $E$, we denote the linear span by $\langle X \rangle$, its closed linear span by $[X]$.

If $I$ is set, $F(I)$ will denote the set of finite subsets of $I$. We turn it into a directed set by inclusion.

\medskip

All tensor products between \cst-algebras in this paper are minimal ones. This implies that the tensor product functionals separate points of the tensor product (and also of its multiplier algebra). The completed tensor products will be denoted by $\ot$. For the tensor product of von Neumann algebras, we use the notation $\ovt$. The flip operator on the tensor product of an algebra with itself will be denoted by $\flip$.

\medskip

The multiplier algebra of a  \cst-algebra $A$ will be denoted by $M(A)$.

Consider two \cst-algebras $A$ and $B$ and a linear map $\rho : A \rightarrow M(B)$. We call $\rho$ strict if it is norm bounded and strictly continuous on bounded sets. If $\rho$ is strict, $\rho$ has a unique linear extension $\bar{\rho} : M(A) \rightarrow M(B)$ which is strictly continuous on bounded sets (see proposition 7.2 of \cite{JK3}). The resulting $\bar{\rho}$ is norm bounded and has the same norm as $\rho$. For $a \in M(A)$, we put $\rho(a) = \bar{\rho}(a)$.

Given two strict linear mappings $\rho : A \rightarrow M(B)$ and $\eta : B \rightarrow M(C)$, we define a new strict linear map $\eta \, \rho : A \rightarrow M(C)$ by  $\eta \, \rho = \bar{\eta} \circ \rho$.
The two basic examples of strict linear mappings are
\begin{itemize}
\item Continuous linear functionals on a \cst-algebra.
\item Non-degenerate $^*$-homomorphism. Recall that a $^*$-homomorphism  $\pi : A \rightarrow M(B)$ is called non-degenerate $\Leftrightarrow$ $ B = [ \, \pi(a)\,b \mid a \in A, b \in B\,]$.
\end{itemize}
All strict linear mappings in this paper will arise as the tensor product of
continuous functionals and/or non-degenerate $^*$-homomorphisms.

\medskip

For $\om \in A^*$ and $a \in M(A)$, we define new elements $a \, \om$ and $\om \, a$ belonging to $A^*$ such that $(a\,\om)(x) = \om(x\,a)$ and
$(\om\,a)(x) = \om(a\,x)$ for $x \in A$.

We also define a  functional $\overline{\om} \in A^*$ such that
$\overline{\om}(x) = \overline{\om(x^*)}$ for all $x \in A$. (Sometimes, $\overline{\om}$ will denote the closure of a densely defined bounded functional, but it will be clear from the context what is precisely meant by $\overline{\om}$).

\medskip

If $A$ and $B$ are \cst-algebras, then the tensor product $M(A) \ot M(B)$ is naturally embedded in $M(A \ot B)$.

\smallskip

We will make extensive use of the leg numbering notation. Let us give an example to illustrate it. Consider three \cst-algebras $A$,$B$ and $C$. Then there exists a unique non-degenerate $^*$-homomorphism $\th : A \ot C \rightarrow M(A \ot B \ot C)$ such that $\th(a \ot c) = a \ot 1 \ot c$ for all $a \in A$ and $c \in C$.

For any element $x \in M(A \ot C)$, we define $x_{13} = \th(x) \in M(A \ot B \ot C)$. It will be clear from the context which \cst-algebra $B$ is under consideration.

If we have another \cst-algebra $D$ and a non-degenerate $^*$-homomorphism
$\de : D \rightarrow M(A \ot C)$, we define the non-degenerate $^*$-homomorphism $\de_{13} : D \rightarrow M(A \ot B \ot C)$ such that
$\de_{13}(d) = \de(d)_{13}$ for all $d \in D$.

\medskip

In this paper, we will also use the notion of a Hilbert \cst-module over a \cst-algebra $A$. For an excellent treatment of Hilbert \cst-modules, we refer to \cite{Lan}.

If $E$ and $F$ are Hilbert \cst-modules over the same \cst-algebra, $\cL(E,F)$ denotes the set of adjointable operators from $E$ into $F$.
When $A$ is a \cst-algebra and $H$ is a Hilbert space, $A \ot H$ will
denote the Hilbert space over $A$, which is a Hilbert \cst-module
over $A$.

For the notion of elements affiliated to a \cst-algebra $A$, we refer to \cite{Baa1}, \cite{Wor6} and \cite{Lan} (these affiliated elements are a generalization of closed densely defined operators in a Hilbert space). For these affiliated elements, there exist notions of self adjointness, positivity and a functional calculus  similar to the notions for closed operators in a Hilbert space. We collected some extra results concerning the functional calculus in \cite{JK6}. Self adjointness will be considered as a part of the definition of positivity. If $\sde$ is a positive element affiliated to a \cst-algebra $A$, $\sde$ is called strictly positive if and only if it has dense range. For such an element $\sde$, functional calculus allows us to define for every $z \in \C$ the power $\sde^z$, which is again affiliated to $A$ (see definition 7.5 of \cite{JK6}).

\medskip

Let $H$ be a Hilbert space. The space of bounded operators on $H$ will be denoted by $B(H)$, the space of compact operators on $H$ by $B_0(H)$. Notice that $M(B_0(H)) = B(H)$.

Let $A$ and $B$ be \cst-algebras and $\pi$ a non-degenerate representation of $A$ on $H$. Consider also $\om \in B_0(H)^*$.

For $a \in M(A)$, we will use the notation $\om(a) := \om(\pi(a)) \in \C$.
For $x \in M(A \ot B)$, we use the notation $(\om \ot \io)(x) := (\om \ot \io)\bigl((\pi \ot \io)(x)\bigr) \in M(B)$.

\medskip

\medskip

Consider a \cst-algebra $A$ and a mapping $\al : \R \rightarrow \text{Aut}(A)$ (where $\text{Aut}(A)$ is the set of $^*$-automorphisms of $A$) such that
\begin{enumerate}
\item $\al_s \, \al_t = \al_{s+t}$ for all $t \in \R$.
\item We have for all $a \in A$ that the function $\R \rightarrow A : t
\rightarrow \al_t(a)$ is norm continuous.
\end{enumerate}
Then we call $\al$ a norm continuous one-parameter group on $A$.
It is then easy to prove that the mapping $\R \rightarrow M(A) : t \mapsto \al_t(a)$ is strictly continuous.

\smallskip

There is a standard way to define for every $z \in \C$ a closed
densely defined linear multiplicative operator $\al_z$ in $A$:
\begin{itemize}
\item The domain of $\al_z$ is by definition the set of elements
$ x \in A$ such that there exists a function $f$ from $S(z)$ into $A$
satisfying
\begin{enumerate}
\item $f$ is continuous on $S(z)$
\item $f$ is analytic on $S(z)^0$
\item We have  that $\al_t(x) = f(t)$ for every $t \in \R$
\end{enumerate}
\item Consider $x$ in the domain of $\al_z$ and $f$ the unique function from $S(z)$ into $A$ such that
\begin{enumerate}
\item $f$ is continuous on $S(z)$
\item $f$ is analytic on $S(z)^0$
\item We have  that $\al_t(x) = f(t)$ for every $t \in \R$
\end{enumerate}
Then we have by definition that $\al_z(x) = f(z)$.
\end{itemize}
where $S(z)$ denotes the strip $\{\,y \in \C \mid \text{Im}\,y \in [0,\text{Im}\,z] \, \}$

\medskip

The mapping $\al_z$ is closable for the strict topology in $M(A)$ and we define the strict closure of $\al_z$ in $M(A)$ by $\overline{\al}_z$. For $a \in D(\overline{\al}_z)$, we put $\al_z(a) := \overline{\al}_z(a)$.

Using the strict topology on $M(A)$, $\overline{\al}_z$ can be constructed from the mapping $\R \rightarrow \text{Aut}(M(A)) : t \rightarrow \overline{\al}_t$ in a similar way as  $\al_z$ is constructed from $\al$.
(See \cite{JK3} or \cite{Ant}, where they used the results in \cite{Zsido} to prove more general results.)

\bigskip

We refer to section 1 of \cite{J-V} for an overview of the necessary weight theory on \cst-algebras. Proper weights are by definition lower semi-continuous weights which are non-zero and
densely defined.  

\smallskip

Let $A$ be a \cst-algebra and $\de : A \rightarrow M(A \ot A)$  a non-degenerate $^*$-homomorphism such that $(\de \ot \io)\de = (\io \ot \de)\de$. Then we call $(A,\de)$ a bi-\cst-algebra. 

Consider a proper weight $\vfi$ on $A$. Then
\begin{itemize}
\item We call $\vfi$ left invariant $\Leftrightarrow$ We have for all $a \in \Mfi^+$ and $\om \in A^*_+$ that $\vfi\bigl((\om \ot \io)\de(a)\bigr) = \om(1) \, \vfi(a)$.
\item We call $\vfi$ right invariant $\Leftrightarrow$ We have for all $a \in \Mfi^+$ and $\om \in A^*_+$ that $\vfi\bigl((\io \ot \om)\de(a)\bigr) = \om(1) \, \vfi(a)$.
\end{itemize}

\smallskip

For some extra information on invariant weights, we refer to section 2 \& 3 of \cite{J-V}.

\bigskip\medskip

\sectie{Reduced locally compact quantum groups}

In this section, we recall the definition of a reduced locally compact qauntum group, as introduced  in \cite{J-V} and list the most important properties of such a reduced quantum group. In a last part, we discuss the reduced dual of such a reduced quantum group. For a detailed exposition, we refer to \cite{J-V}.

\medskip

Let us first start with the definition of a \cst-algebraic quantum group.

\begin{definition} \label{red.def1}
Consider a \cst-algebra $A$ and a non-degenerate $^*$-homomorphism $\de : A \rightarrow M(A \ot A)$ such that 
\begin{itemize}
\item $(\de \ot \io)\de = (\io \ot \de)\de$.
\item $A = [\,(\om \ot \io)\de(a) \mid \om \in A^*,a \in A \, ] = [\,(\io \ot \om)\de(a) \mid \om \in A^*, a \in A \, ]$.
\end{itemize}
Assume moreover the existence of
\begin{itemize}
\item A faithful left invariant approximate KMS weight $\vfi$ on
$(A,\de)$.
\item A right invariant approximate KMS weight $\psi$ on $(A,\de)$.
\end{itemize}
Then we call $(A,\de)$ a reduced \cst-algebraic quantum group.
\end{definition}
 
The weak KMS property is a weaker condition than the usual KMS property for a weight on a \cst-algebra but it turns out that every proper left or right invariant weight on such a reduced \cst-algebraic quantum group is automatically faithful and KMS. Moreover, proper left invariant weights are unique up to a scalar (and similarly for proper right invariant weights).

\bigskip\medskip

For the rest of this paper, we will fix a reduced \cst-algebraic quantum group $(A,\de)$ together with a faithful left invariant KMS weight $\vfi$ on $(A,\de)$ such that there exists a GNS-construction $(H,\io,\la)$ for $\vfi$ (here $\io$ denotes the identity map of $A$).

So we assume (for convenience purposes) that $A$ acts on the GNS-space of its left Haar weight $\vfi$ in a particular way. This is obviously not very essential. We let $\tilde{A}$ denote the von Neumann algebra acting on $H$ generated by $A$.

\medskip

Let us give a short overview of the main objects associated to our \cst-algebraic quantum group $(A,\de)$:
\begin{trivlist}
\item[\ \bf 1.]  \bf The antipode and its polar decomposition: \rm

\medskip

The antipode $S$ of $(A,\de)$ is a closed linear mapping in $A$ determined by the following properties:
\begin{itemize}
\item We have for all $a,b \in \Nfi$ that
$$(\io \ot \vfi)(\de(a^*)(1 \ot b)) \in D(S)$$
and
$$S\bigl((\io \ot \vfi)(\de(a^*)(1 \ot b))\bigr) = (\io \ot \vfi)((1 \ot a^*)\de(b)) \ .$$
\item The set 
$$\langle \, (\io \ot \vfi)(\de(a^*)(1 \ot b))\mid a,b \in \Nfi \, \rangle$$
is a core for $S$.
\end{itemize}

There exists a unique $^*$-antiautomorphism $R$ on $A$ and a unique norm continuous one-parameter group $\tau$ on $A$ such that
\begin{itemize}
\item $R^2=\io$,
\item $R$ and $\tau$ commute,
\item $S = R \, \tau_{-\frac{i}{2}}$.
\end{itemize}
The pair $R,\tau$ is called the polar decomposition of $S$. The $^*$-antiautomorphism $R$ is called the unitary antipode of $(A,\de)$ and the one-parameter group $\tau$ is called the scaling group of $(A,\de)$.

\medskip\medskip

\item[\ \bf 2.]  \bf The Haar weights and their modular groups:\rm

\medskip

The unitary antipode $R$ satisfies the equality $\flip(R \ot R)\de$. So we can define the right invariant faithful KMS weight $\psi$ on $(A,\de)$ as $\psi = \vfi R$. 

The modular group of $\vfi$ is denoted by $\si$, the modular group of $\psi$ is denoted by $\sip$. These one-parameter groups are related by the formula $\sip_t = R \si_{-t} R$ for $t \in \R$.

\medskip

The different relations between $\si$, $\sip$ and $\tau$ are collected in the following list.
\begin{itemize}
\item The automorphism groups $\si,\sip$ and $\tau$ commute pairwise.
\item We have the following commutation relations for all $t \in \R$:
$$\begin{array}{rclcrcl}
\de \, \si_t  & =  & (\tau_t \ot \si_t)\de  & \hspace{2cm}  & \de\, \sip_t & = & (\sip_t \ot \tau_{-t}) \de \\
\de \,\tau_t & =  &(\tau_t \ot \tau_t) \de & \hspace{2cm}  & \de \, \tau_t & = & (\si_t \ot \sip_{-t}) \de
\end{array}$$
\item There exists a  number $\nu > 0$ such that
$$\begin{array}{rclcrcl}
\vfi \, \sip_t & = & \nu^t \, \vfi &\hspace{1.5cm}  & \psi \, \si_t  &=  &\nu^{-t} \,\psi  \\
 \psi\, \tau_t & = & \nu^{-t} \, \psi & \hspace{1.5cm} & \vfi \, \tau_t  &=  &\nu^{-t} \, \vfi
\end{array}$$
for all $t \in \R$.
\end{itemize}

The number $\nu$ is called the scaling constant of $(A,\de)$. It is not clear yet whether this number can be different from 1.

\medskip\medskip

\item[\ \bf 3.]  \bf The modular element:\rm

\medskip

In the next part, we will use the terminology and notations of section 1.4 of \cite{J-V}. There exists a unique strictly positive element $\sde$ affiliated to $A$ such that  
$\si_t(\sde) = \sip_t(\sde) = \nu^t \, \sde$ for all $t \in \R$ and $\psi = \vfi_\sde$.
So $\sip_t(x) = \sde^{it} \, \si_t(x) \, \sde^{-it}$ for all $t \in \R$ and $x \in A$.

\smallskip

We use the equality $\psi = \vfi_\sde$ to define a GNS-construction $(H,\io,\ga)$ for $\psi$ such that $\ga = \la_{\sde}$.

\medskip

Let us list some elementary properties of $\sde$:
\begin{itemize}
\item $\de(\sde) = \sde \ot \sde$
\item $\tau_t(\sde) =  \sde$  for $t \in \R$ and $R(\sde) = \sde^{-1}$.
\item Let $t \in \R$. Then $\sde^{it}$ belongs to $D(\bar{S})$ and
$S(\sde^{it}) = \sde^{-it}$.
\end{itemize}

\medskip\medskip

\item[\ \bf 4.] \bf The multiplicative unitary: \rm

\medskip

The multiplicative unitary $W$ of $(A,\de)$ (in this particular GNS-construction $(H,\io,\la)$)
is the unitary element in $B(H \ot H)$ such that $W(\la \ot \la)(\de(b)(a \ot 1)) = \la(a) \ot \la(b)$ for all $a,b \in \Nfi$.

The operator $W$ satisfies the Pentagonal equation: $W_{12} W_{13} W_{23} = W_{23} W_{12}$. Moreover, it encodes all the information about $(A,\de)$ in the following way:
\begin{itemize} 
\item $A = [\,(\io \ot \om)(W) \mid \om \in B_0(H)^* \, ]$, 
\item $\de(x) = W^* (1 \ot x)W$ for all $x \in A$.
\end{itemize} 
\end{trivlist}

\medskip

The main aim of this paper is to show that the \lq universal\rq\ quantum group which is associated to $(A,\de)$ has this same rich analytical structure.

\bigskip\bigskip

Given such a reduced quantum group $(A,\de)$, there is a standard way to construct the dual $(\ah,\deh)$ of $(A,\de)$. The pair $(\ah,\deh)$ is again a reduced \cst-algebraic quantum group
and can be easily defined in terms of the multiplicative unitary:
\begin{itemize}
\item $\ah = [ \, (\om \ot \io)(W) \mid \om \in B_0(H)^* \,]$,
\item $\deh(x) =  \Sigma W(x \ot 1)W^* \Sigma$ for all $x \in \ah$,
\end{itemize}
where $\Sigma$ denotes the flip map on $H \ot H$. The multiplicative unitary $W$ belongs to $M(A \ot \ah)$. 

\smallskip

The symbols for the  objects associated to the quantum group $(\ah,\deh)$ (antipode, unitary antipode,...) will be obtained by adding $\ \hat{}\ $ to the symbol of the counterpart on the level of $(A,\de)$ (e.g. the scaling group of $(\ah,\deh)$ will be denoted by $\tauh$).

\bigskip\medskip

Using the multiplicative unitary to define the dual somewhat hides the fact that the dual can be obtained  from a construction which resembles the construction of the reduced group \cst-algebra of a locally compact group. Let us strengthen the analogy with the group case by introducing  the closed subspace $\dual$ of $A^*$:
$$\dual = [ \, a \vfi b^* \mid a,b \in \Nfi\, ] = [\,\om\restr_A \mid \om \in B(H)_*\,]\ .$$

The topological dual $A^*$ is a Banach algebra under the multiplication $A^* \times A^* \rightarrow A^* : (\om,\th) \mapsto \om\,\th$ given by $(\om\, \th)(x) = (\om \ot \th)\de(x)$ for all $\om,\th \in A^*$  and $x \in A$. The set $\dual$ is a two sided ideal in $A^*$.

We define the injective contractive algebra homomorphism $\lambda : A^* \rightarrow M(\ah)$ such that $\lambda(\om) = (\om \ot \io)(W)$ for $\om \in \dual$. Then $\lambda(\dual)$ is a dense 
subalgebra of $\ah$.

\medskip

If $S$ is unbounded, the algebra $\dual$ does not carry an appropriate $^*$-structure. It is however possible to find an subalgebra of $\dual$ which carries a $^*$-structure:

\smallskip

Define the subspace $\duals$ of $\dual$ as 
$$\duals = \{\, \om \in \dual \mid \exists \,\th \in \dual : \th(x) = \overline{\om}(S(x)) \text{ for all } x \in D(S) \, \}\ .$$
We define the antilinear mapping $.^* : \duals \rightarrow \duals$ such that
$\om^*(x) = \overline{\om}(S(x))$  for all $\om \in \duals$ and $x \in D(S)$.
Then $\duals$ is a subalgebra of $\dual$ and becomes a $^*$-algebra under the operation $.^*$\ 

\bigskip\medskip

\sectie{The generator of the universal representation of the dual}

In the last part of the previous section, we introduced the $^*$-algebra $\duals$ as a sub algebra of $L^1(A)$. On this $^*$-algebra, we introduce the natural norm $\|.\|_*$ such that 
$$\|\om\|_* = \max\{\|\om\|,\|\om^*\|\}$$
for all $\om \in \duals$ \ (where $\|.\|$ denotes the norm on $A^*$). It is then easy to check that $\duals$ together with this norm $\|.\|_*$ becomes a Banach $^*$-algebra. Whenever we use topological concepts connected to $\duals$ without further mention, we will always be working with the norm $\|.\|_*$. 

\medskip

As a consequence, we can form the universal enveloping \cst-algebra $\ahu$ of $\duals$.  Recall that $\ahu$ is formed in the following way. First one defines a norm $\|.\|_{\text{\tiny u}}$ on $\duals$ such that 
$$\|\om\|_{\text{\tiny u}} =  \sup \{\,\|\th(\om)\| \mid \th \text{ a  $^*$-representation of } \duals \text{ on a Hilbert space} \, \}$$ 
for all $\om \in \duals$ \ (because $\lambda$ is an injective $^*$-representation, one gets a norm and not merely a semi-norm). 

\medskip

In a next step, one defines $\ahu$ to be a completion of $\duals$ with respect to this norm $\|.\|_{\text{\tiny u}}$. The embedding of $\duals$ into $\ahu$ will be denoted by $\lambdau$. Then the pair $(\ahu,\lambdau)$ is (up to a $^*$-isomorphism) determined by the following universal property:

\smallskip

Let $C$ be any \cst-algebra and $\th : \duals \rightarrow C$ a $^*$-representation. Then there exists a unique $^*$-homomorphism $\th_{\text{\tiny u}} : \ahu \rightarrow C$ such that $\th_{\text{\tiny u}} \, \lambdau = \th$.

\medskip

By  choosing the completion $\ahu$ in the right way, we can assume that $\ahu$ acts on  a Hilbert space $H_{\text{\tiny u}}$.

\bigskip\medskip

In this section, we will prove the existence of a unitary  element $\hat{\cV} \in M(A \ot \ahu)$ such that $\lambdau(\om) = (\om \ot \io)(\hat{\cV})$ for all $\om \in \duals$ \ ($\hat{\cV}$ is called the generator for $\lambdau$). This will immediately imply that a similar property holds for any $^*$-representation of $\duals$. 

This result was proven by Kirchberg for Kac algebras and a careful analysis of his proofs shows that they can be easily transformed to proofs of the result in the general quantum group case. In the rest of this section we will give the transformed proofs and indicate what had to be changed to them. We will essentially follow  the discussion in sections 1.4 and 3.1 of \cite{E}.

\bigskip

The main \lq problem\rq\ in transforming the proofs from the Kac algebra setting to the general quantum group setting stems from the following fact.

\medskip

In the Kac algebra framework, we have that $\duals = \dual$ as Banach spaces and this is not the case in the general quantum group setting. But we know that $\dual$ is isomorphic to the predual $\tilde{A}_*$ implying that $\dual^*$ is isomorphic to $\tilde{A}$.  In section 1.4 of \cite{E}, the product in $\tilde{A}$ is then used to define the Kronecker product between two $^*$-representations of $\duals$.

If $\th \in \ahu^*$, we know in general a priori only that $\th \lambdau$ in $\duals^*$, it is not clear that $\th \lambdau $ can be extended to an element in $\dual^*$ and in this way give rise to an element in $\tilde{A}$ \ (once we have the generator $\hat{\cV}$ at our disposal, this is obvious). 

\medskip

But it will turn out that in order to define the Kronecker product $\lambda \times \lambdau$ (which will be sufficient to prove the existence of the generator), it is enough   to define a module action of a well-behaved subset of $A$ on $\duals^*$. 

\begin{lemma} \label{gen.lem1}
Consider $\om \in \duals$ and $x \in D(S)$. Then 
\begin{enumerate}
\item $\om x$ and $x \om$ belong to $\duals$.
\item $(\om x)^* = \om^* S(x)^*$ and $(x \om)^* = S(x)^* \om$.
\item $\|\om x\|_*,\|x \om \|_* \leq \max\{\|x\|,\|S(x)\|\}\,\,\|\om\|_*$.
\end{enumerate}
\end{lemma}
\begin{demo}
By definition of $\om^*$ and using proposition 5.22 of \cite{J-V}, we have for all $y \in D(S)$ that 
$$\overline{(\om x)(S(y)^*)} = \overline{\om(x S(y)^*)}
= \overline{\om(S(S(x)^* y)^*)}  = \om^*(S(x)^* y) = (\om^* S(x)^*)(y)\ .$$
By definition of $\duals$ and its $^*$-operation, this implies that $\om x$ belongs to $\duals$ and $(\om x)^* = \om^* S(x)^*$. Then we have also immediately that $\|\om x\|_*  \leq \max\{\|x\|,\|S(x)\|\}\,\,\|\om\|_*$. 
The result about $x \om$ is proven in a similar way.
\end{demo}

\begin{remark} \rm
This lemma  implies that we can define the following module operations on $\duals^*$. 

Consider $F \in \duals^*$ and $x \in D(S)$. Then we define $F x$,$x F$ $\in \duals^*$ such that $(F x)(\om) = F(\om x)$ and $(x  F)(\om) = F(x \om)$ for all $\om \in \duals^*$. Of course, $\|F x\|,\|x F\| \leq \max\{\|x\|,\|S(x)\|\}\,\,\|F\|$.

In this way, $\duals^*$ becomes a bimodule over $D(S)$. Although we will not use this module notations anymore, it is implicitly present in proposition \ref{gen.prop1}
\end{remark}

\bigskip\bigskip

Let us introduce a symbol to denote the pullback of the map $\lambda : \dual \rightarrow \ah : \om \mapsto (\om \ot \io)(W)$.

\begin{notation}
We define the linear contraction $\lambda^* : B_0(H)^* \rightarrow A$ such that
$\lambda^*(\th)= (\io \ot \th)(W)$ for all $\th \in B_0(H)^*$.
\end{notation}

\begin{remark} \label{gen.rem1} \rm
Notice that $\om(\lambda^*(\th)) = \th(\lambda(\om))$ for all $\om \in \duals$ and $\th \in B_0(H)^*$ so that $\lambda^*$ is really the pullback of $\lambda$. By proposition 8.3 of \cite{J-V}, we get for all $\th \in B_0(H)^*$ that $\lambda^*(\th) \in D(S)$ and $S(\lambda^*(\th))^* = \lambda(\bar{\th})$.

Consider $\om \in \duals$. By lemma \ref{gen.lem1}, the following properties hold.
\begin{itemize}
\item We have for all $\eta \in B_0(H)^*$ that $\om\,\lambda^*(\eta) \in \duals$
and $(\om\,\lambda^*(\eta))^* = \om\,\lambda^*(\bar{\eta})$.
\item The linear map $B_0(H)^* \rightarrow \duals : \eta \mapsto \om\,\lambda^*(\bar{\eta})$ is continuous.
\end{itemize}
\end{remark}

\bigskip

We want to mimic the proofs of proposition 1.4.2 and theorem 1.4.3 of \cite{E} (which are due to Kirchberg) to define a new $^*$-representation $\mu$ of $\duals$ on $B(H \ot H_{\text{\tiny u}})$. The proof of theorem 1.4.2 requires the product in $\tilde{A}$ but it turns out that in this case the module action of $D(S)$ on $\duals^*$ is sufficient. The map $\mu$ in the next proposition is nothing els but the Kronecker product $\lambda \times \lambdau$.

\begin{proposition} \label{gen.prop1}
There exists a unique $^*$-representation $\mu : \duals \rightarrow B(H \ot H_{\text{\tiny u}})$
such that 
$$\langle \mu(\om)(v_1 \ot w_1), v_2 \ot w_2 \rangle
=  \langle \lambdau\bigl(\om\,\lambda^*(\om_{v_1,v_2})\bigr)\,w_1 , w_2 \rangle \ .$$
for all $v_1,v_2 \in H$ and $w_1,w_2 \in H_{\text{\tiny u}}$.
\end{proposition}
\begin{demo}
Fix an orthonormal basis $(e_k)_{k \in K}$ for $H$. For every $k,l \in K$, we define $x_{kl} = \lambda^*(\om_{e_l,e_k})$. 

\smallskip

Take a finite subset $L$ of $K$ and for every $l \in L$ a vector $w_l \in H_{\text{\tiny u}}$.

\medskip

Choose $\om \in \duals$. Then we have that
\begin{eqnarray*}
\sum_{k \in K} \ \left\|\,\sum_{l \in L} \lambdau(\om x_{kl}) w_l \,\right\|^2  & = &  \sum_{k \in K} \,\, \sum_{l,l' \in L} \langle \lambdau(\om x_{kl}) w_l , 
\lambdau(\om x_{k l'}) w_{l'} \rangle \\
& = & \sum_{k \in K} \,\, \sum_{l,l' \in L} \langle \lambdau\bigl((\om x_{k l'})^*(\om x_{k l})\bigr) w_l , w_{l'} \rangle \ .
\end{eqnarray*}
Using  remark \ref{gen.rem1} , this implies that
\begin{equation}
\sum_{k \in K}\ \left\|\,\sum_{l \in L} \lambdau(\om x_{kl}) w_l \,\right\|^2
= \sum_{k \in K} \,\, \sum_{l,l' \in L} \langle \lambdau\bigl((\om^* x_{l'k})(\om x_{kl})\bigr) w_l , w_{l'} \rangle \ . \label{gen.eq1}
\end{equation}

\medskip

Fix $l,l' \in L$ for the moment. 
We have for every set $M \in F(K)$ 
$$\sum_{k \in M} x_{l'k} \ot x_{kl} = \sum_{k \in M} (\io \ot \io \ot \om_{e_k,e_{l'}})(W_{13})\, (\io \ot \io \ot \om_{e_l,e_k})(W_{23})\ .$$
Hence lemma 9.5 of \cite{J-V} implies that the net $\bigl(\,\sum_{k \in M} x_{l'k} \ot x_{kl}\,\bigr)_{M \in F(K)}$ is bounded and converges strictly  to  $(\io \ot \io \ot \om_{e_l,e_{l'}})(W_{13} W_{23})$ in $M(A \ot A)$. But the Pentagonal equation implies that this last expression is equal to 
$$(\io \ot \io \ot \om_{e_l,e_{l'}})(W_{12}^* W_{23} W_{12})
= W^* (1 \ot (\io \ot \om_{e_l,e_{l'}})(W))\, W = \de(x_{l'l}) \ .$$
So we conclude that the net $\bigl(\,\sum_{k \in M} x_{l'k} \ot x_{kl}\,\bigr)_{M \in F(K)}$ is bounded and converges strictly to $\de(x_{l'l})$.

Since 
$$\sum_{k \in M} (\om^* x_{l'k})(\om x_{kl}) 
= \bigl(\,(\om^* \ot \om)\bigl[\,\sum_{k \in M} x_{l'k} \ot x_{kl}\,\bigr]\,\bigr)\,\de$$
for all $M \in F(K)$, we conclude that the net $\bigl(\,\sum_{k \in M} (\om^* x_{l'k})(\om x_{kl})\,\bigr)_{M \in F(K)}$ converges in $\dual$ to 
$((\om^* \ot \om)[\de(x_{l'l})])\,\de = (\om^*  \om)\,x_{l'l}$.

\medskip

But we have for every $M \in F(K)$ also that
$$\bigl(\,\sum_{k \in M} (\om^* x_{l'k})(\om x_{kl}) \,\bigr)^*
= \sum_{k \in M} (\om^* x_{l k})(\om x_{kl'}) \ , $$
implying that the net $\bigl(\,\bigl(\,\sum_{k \in M} (\om^* x_{l'k})(\om x_{kl})\,\bigr)^* \,\bigr)_{M \in F(K)}$ converges in $\dual$ to 
$(\om^* \om)\,x_{l l'}$, which is equal to $((\om^* \om)\,x_{l'l})^*$.
So we conclude that the net $\bigl(\,\sum_{k \in M} (\om^* x_{l'k})(\om x_{kl})\,\bigr)_{M \in F(K)}$ converges in $\duals$ to $(\om^*  \om)\,x_{l'l}$.

\medskip

Because $\lambdau : \duals \rightarrow \ahu$ is bounded, we now conclude from equation  \ref{gen.eq1} that
\begin{equation} \label{gen.eq2}
\sum_{k \in K}\ \left\|\,\sum_{l \in L} \lambdau(\om x_{kl}) w_l \,\right\|^2
= \sum_{l,l' \in L} \langle \lambdau((\om^* \om) \,x_{l'l}) w_l , w_{l'} \rangle \ . 
\end{equation}

\bigskip

Define the linear functional $\eta$ on $\duals$ such that 
$$\eta(\om) = \sum_{l,l' \in L} \langle \lambdau(\om\,x_{l'l}) w_l , w_{l'} \rangle$$ 
for all $\om \in \duals$. Then equation \ref{gen.eq2} implies that $\eta$ is a positive functional on $\duals$.

\medskip

Take $\om \in \duals$. Then 
\begin{eqnarray*}
|\eta(\om)|^2 & = & \left|\ \sum_{l,l' \in L} \langle \lambdau(\om \,x_{l'l}) w_l , w_{l'} \rangle\ \right|^2
\\ & = & \left|\ \sum_{l' \in L }\, \langle\, \sum_{l \in L}  \lambdau(\om \,x_{l'l})  w_l , w_{l'} \rangle\ \right|^2 \\
& \leq &  \left(\ \sum_{l' \in L} \, \left\|\,\sum_{l \in L}  \lambdau(\om \,x_{l'l}) \, w_l\,\right\|^2\ \right)
 \ \bigl(\,\sum_{l' \in L} \|w_{l'}\|^2\,\bigr) 
\end{eqnarray*}
where in the last inequality, we used the Cauchy Schwarz inequality in $\oplus_{l' \in L} \, H_{\text{\tiny u}}$. Hence equation \ref{gen.eq2} implies that
\begin{eqnarray*}
|\eta(\om)|^2 & \leq & \left(\ \sum_{l' \in K} \, \left\|\,\sum_{l \in L}  \lambdau(\om \,x_{l'l}) \, w_l\,\right\|^2\ \right)  \ \bigl(\,\sum_{l' \in L} \|w_{l'}\|^2\,\bigr) 
 \\
& = & \bigl(\,\sum_{l,l' \in L} \,\langle \lambdau((\om^* \om) \,x_{l'l}) w_l , w_l' \rangle\,\,\bigr)
\ \bigl(\,\sum_{l' \in L} \|w_{l'}\|^2 \, \bigr) \\
& = & \bigl(\,\sum_{l' \in L} \|w_{l'}\|^2 \, \bigr) \ \eta(\om^* \om) \ .
\end{eqnarray*}
Therefore theorem 37.11 of \cite{Bon} implies that $\eta$ is continuous and has norm less than $\sum_{l \in L} \|w_l\|^2$. So we conclude from equation \ref{gen.eq2} 
$$ \sum_{k \in K}\ \left\|\,\sum_{l \in L} \lambdau(\om \,x_{kl}) w_l \,\right\|^2 = \eta(\om^* \om)  \leq  \bigl(\,\sum_{l \in L} \|w_l\|^2\,\bigr) \ \|\om^* \om\| \leq \left\|\,\sum_{l \in L} e_l \ot w_l\,\right\|\ \|\om\|^2\ .$$

\medskip

From this all, we get the existence of a contractive linear map $\mu : \duals \rightarrow B(H \ot H_{\text{\tiny u}})$ such that
$$\mu(\om)(e_l \ot w) = \sum_{k \in K} e_k \ot \lambdau(\om \, x_{kl}) \, w$$
for all $l \in K$ and $w \in H_{\text{\tiny u}}$.

\medskip

So we get for $k,l \in K$ and $w_1,w_2 \in H$ that
$$\langle \mu(\om)(e_l \ot w_1) , e_k \ot w_2 \rangle = \langle \lambdau\bigl(\om \,\lambda^*(\om_{e_l,e_k})\bigr) w_1 , w_2 \rangle \ .$$
Since the linear function $B_0(H)^* \rightarrow \duals : \th \rightarrow \om \lambda^*(\th)$ is continuous (see remark \ref{gen.rem1}), we conclude  that
$$\langle \mu(\om)(v_1 \ot w_1), v_2 \ot w_2 \rangle
=  \langle \lambdau\bigl(\om\,\lambda^*(\om_{v_1,v_2})\bigr)\,w_1 , w_2 \rangle \ .$$
for all $v_1,v_2 \in H$ and $w_1,w_2 \in H_{\text{\tiny u}}$.

\medskip

In the last part of this proof, we show that $\mu$ is a $^*$-homomorphism.
\begin{enumerate}
\item $\mu$ is selfadjoint:

Take $\om \in \duals$. Choose $v_1,v_2 \in H$ and $w_1,w_2 \in H_{\text{\tiny u}}$. Then 
\begin{eqnarray*}
& & \langle \mu(\om)^* (v_1 \ot w_1), v_2 \ot w_2 \rangle
= \langle  v_1 \ot w_1, \mu(\om)(v_2 \ot w_2) \rangle \\
& & \spat = \overline{\langle \mu(\om)(v_2 \ot w_2) , v_1 \ot w_1 \rangle}
= \overline{\langle \lambdau(\om\,\lambda^*(\om_{v_2,v_1}))\,w_2 , w_1 \rangle} \\
& & \spat = \langle w_1 , \lambdau\bigl(\om\,\lambda^*(\om_{v_2,v_1})\bigr)\,w_2 \rangle
= \langle  \lambdau\bigl(\om\,\lambda^*(\om_{v_2,v_1})\bigr)^* w_1 , w_2 \rangle \\
& & \spat = \langle \lambdau\bigl((\om\,\lambda^*(\om_{v_2,v_1}))^*\bigr) w_1 , w_2 \rangle
= \langle \lambdau\bigl(\om^* \, \lambda^*(\om_{v_1,v_2})\bigr) w_1 , w_2 \rangle \\
& & \spat = \langle \mu(\om^*) (v_1 \ot w_1) , v_2 \ot w_2 \rangle \ .
\end{eqnarray*}
So we conclude that $\mu(\om)^*  = \mu(\om^*)$.

\item $\mu$ is multiplicative:

Choose $\om \in \duals$. By equation \ref{gen.eq2}, we have for every finite subset $L$ of $K$ and vectors \newline $w_l \in H_{\text{\tiny u}}\ (l \in L)$ that 
$$\langle \mu(\om) \,\bigl(\,\sum_{l \in L} e_l \ot w_l\,\bigr)  , \mu(\om) \,\bigl(\,\sum_{l \in L} e_l \ot w_l\,\bigr) \rangle
=  \langle \mu(\om^* \om)\, \bigl(\,\sum_{l \in L} e_l \ot w_l\,\bigr) , \sum_{l \in L} e_l \ot w_l \rangle \ .$$
Hence $\mu(\om^* \om) = \mu(\om)^* \mu(\om)
= \mu(\om^*) \, \mu(\om)$,
where we used the selfadjointness of $\mu$ in the last equality. Polarization yields that 
$\mu(\th^* \om) = \mu(\th^*) \mu(\om)$ for all $\om,\th \in \duals$.
\end{enumerate}
\end{demo}   
 
\begin{remark}\rm \label{gen.rem2}
In the notation of the previous proposition, we get for all $v,w \in H$ and $\om \in \duals$ that  $$(\om_{v,w} \ot \io)(\mu(\om)) = \lambdau\bigl(\om\,\lambda^*(\om_{v,w}\bigr)) \in \ahu \ .$$
\end{remark}

\bigskip

For every $\om \in \dual$, we will denote the unique normal functional in $\tilde{A}_*$ which extends $\om$ by $\tilde{\om}$. As discussed before, it might still be possible that there exists an element $\th \in \ahu^*$ such that there does not exist an element $y \in \tilde{A}$ satisfying $\th(\lambdau(\om)) = \tilde{\om}(y)$ for all $\om \in \duals$. 

For this reason, we have to adapt the proofs of \cite{E} a little bit further. In the next lemma, we will provide sufficiently many elements $\th$ for which there do exist such elements $y$ as mentioned above.

\medskip\medskip

Let $t \in \R$. Since $S$ and $\tau_t$ commute, we have for every  $\om \in \duals$ that the element $\om \,\tau_t$ belongs to $\duals$ and $(\om \, \tau_t)^* = \om^* \, \tau_t$.  \inlabel{gen.eq7}
 
Combining this  with the fact that $(\tau_t \ot \tau_t)\de = \de$, we see that the mapping $\duals \rightarrow \duals : \om \rightarrow \om \tau_t$ is a $^*$-automorphism. 

\smallskip\smallskip

Also notice that equation \ref{gen.eq7} implies that the mapping $\R \rightarrow L^1_*(A) : t \mapsto \om \, \tau_t$ is continuous (see also proposition 8.23 of \cite{J-V}).

\medskip\medskip

Due to the universal property of $\ahu$, we can therefore introduce the following norm continuous one-parameter group on $\ahu$. It will be the scaling group of the universal dual of $(A,\de)$.

\begin{definition}
There exists a unique norm continuous one-parameter group $\taudu$ on $\ahu$ such that 
\newline $\taudu_t(\lambdau(\om)) = \lambdau(\om \tau_{-t})$ for all $t \in \R$ and $\om \in \duals$.
\end{definition}

Notice that norm continuity follows because $\lambdau(\duals)$ is dense in $\ahu$, the map $\lambdau : \duals \rightarrow \ahu$ is a contraction and the function
$\R \rightarrow \duals : t \rightarrow \om \tau_t$ is norm continuous for every $\om \in \duals$

\begin{lemma}
Consider $\th \in \ahu^*$ and $n \in \N$. Define $\th_n \in \ahu^*$ such that 
$$\th_n(x) = \frac{n}{\sqrt{\pi}} \int \exp(-n^2 t^2) \, \th(\taudu_t(x)) \, dt$$
for all $x \in \ahu$.
Then there exists a unique element $y \in \tilde{A}$ such that 
$\th_n(\lambdau(\om)) = \tilde{\om}(y)$ for all $\om \in \duals$.
\end{lemma}
\begin{demo}
We have for all $\om \in \duals$ that 
\begin{eqnarray*}
\th_n(\lambdau(\om)) & = & \frac{n}{\sqrt{\pi}} \int \exp(-n^2 t^2) \, \th\bigl(\taudu_t(\lambdau(\om))\bigr) \, dt \\
& = & \frac{n}{\sqrt{\pi}} \int \exp(-n^2 t^2) \, \th(\lambdau(\om\tau_{-t})) \, dt \\
& = & \th\left(\ \frac{n}{\sqrt{\pi}} \int \exp(-n^2 t^2)\,\om\,\tau_{-t}\, dt\ \right) \ .
\end{eqnarray*}
Define the function $F : \duals \rightarrow \duals : \om \mapsto \frac{n}{\sqrt{\pi}} \int \exp(-n^2 t^2)\,\om\,\tau_{-t}\, dt$. Then this function is clearly $\|.\|,\!\|.\|$ continuous.
We have moreover for all $\om \in \duals$ that (see e.g. the proof of lemma 8.33 of \cite{J-V})
$$F(\om)^* =   \frac{n}{\sqrt{\pi}} \int \exp(-n^2 (t-\frac{i}{2})^2)\,\,\overline{\om} R \tau_{-t}\, dt\ ,$$
which implies that the function $\duals \rightarrow \duals : \om \mapsto F(\om)^*$ is also 
$\|.\|,\!\|.\|$ continuous. Consequently, the function $F$ is $\|.\|,\!\|.\|_*$ continuous.

\smallskip

But we have for all $\om \in \duals$ that $\th_n(\lambdau(\om)) = \th\bigl(\lambdau(F(\om))\bigr)$, implying  the existence of a unique element $\eta \in \dual^*$ such that $\eta(\om) = \th_n(\lambdau(\om))$ for all $\om \in \duals$. Using the natural isomorphisms $L_1(A)^* \cong (\tilde{A}_*)^* \cong \tilde{A}$,  there exists  an element $y \in \tilde{A}$ satisfying $\eta(\om) = \tilde{\om}(y)$ for all $\om \in \dual$ and the lemma follows.
\end{demo}

\begin{corollary} \label{gen.cor1}
The set
$$\{\, \th \in \ahu^* \mid \exists y \in \tilde{A}, \forall \om \in \duals : \th(\lambdau(\om)) = \tilde{\om}(y) \,\}$$
is separating for $\ahu$.
\end{corollary}

Notice that in order for things to work in the previous discussion , we have to stick with the \cst-algebra $\au$ and not go to the universal enveloping von Neumann algebra of $\au$ as is done in \cite{E} (to be more precise, the one-parameter group $\taudu$ can be point wisely extended to the enveloping von Neumann algebra but this extension does not satisfy any obvious continuity property).  Remark \ref{gen.rem2} will allow us to stick to the \cst-algebra setting.

\medskip

Using the universal property of $(\ahu,\lambdau)$, we define $^*$-homomorphisms $s_\lambda : \ahu \rightarrow \ah$ and $s_\mu : \ahu \rightarrow B(H \ot H_{\text{\tiny u}})$  such that 
$s_\lambda \, \lambdau = \lambda$ and $s_\mu \,\lambdau = \mu$.

\smallskip

Thanks to corollary \ref{gen.cor1}, we can use the techniques in the proof of lemma 3.1.1 and proposition 3.1.3 of \cite{E} to get closer to our goal. First we need some extra results from section 8 of \cite{J-V}.

\smallskip

We define the subspace $\cI$ of $\dual$ as follows (notation 8.4 of \cite{J-V}):
$$\cI =  \{ \, \om \in \dual \mid \text{There exists a number } M \geq 0 \text{\ s.t.\ } |\om(x^*)| \leq M \, \|\la(x)\| \text{ for all } x \in \Nfi \, \}\ .$$

By Riesz' theorem for Hilbert spaces, there exists for every $\om \in \cI$ a unique element $\xi(\om) \in H$ such that $\om(x^*) = \langle \xi(\om) , \la(x) \rangle$ for $x \in \Nfi$.
Result 8.6 of \cite{J-V} tells us that $\cI$ is a left ideal in $\dual$ and that $\xi(\th \om) = \lambda(\th) \, \xi(\om)$ for all $\th \in \dual$ and $\om \in \cI$. 

\medskip

\begin{lemma}
We have that $\ker s_{\lambda} \subseteq \ker s_\mu$.
\end{lemma}
\begin{demo}
Choose $p \in \ker s_{\lambda}$. 

Take $\eta \in \cI \cap \duals$. Choose also $a,b \in \Nfi$, $c \in \Npsi$ and define $v,w \in H$ by $v = J \la(c^* a )$, $w = J \la(b)$. Since $\hat{R}$ is implemented by $J$ (see proposition 8.17 of \cite{J-V}) and $(R \ot \hat{R})(W) =  W$ (see the remarks before proposition 8.18 of \cite{J-V}), we get 
\begin{eqnarray*}
& & \lambda^*(\om_{v,w})  =  (\io \ot \om_{v,w})(W) = R\bigl((\io \ot \om_{J w , J v})(W)\bigr)
\\ & & \spat = R\bigl((\io \ot \om_{\la(b ),\la(c^* a)})(W)\bigr) 
 =  R\bigl((\io\ \ot \vfi)(\de(a^* c)(1 \ot b))\bigr) \ .
\end{eqnarray*}
By the right invariant version of result 2.6 of \cite{J-V}, $(\io\ \ot \vfi)(\de(a^* c)(1 \ot b))$ belongs to $\Npsi$. Therefore the previous equation and the fact that $\psi = \vfi R$ imply that
$\lambda^*(\om_{v,w})$ belongs to $\Nfi^*$.

\smallskip

By remark \ref{gen.rem2}, we have for all  $\om \in \duals$ that 
$$(\om_{v,w} \ot \io)(s_\mu(\lambdau(\om))\, \mu(\eta)) = (\om_{v,w} \ot \io)(\mu(\om) \mu(\eta))
= (\om_{v,w} \ot \io)(\mu(\om \eta)) \in \ahu \ .$$
As a consequence, we find that $(\om_{v,w} \ot \io)(s_\mu(x) \,\mu(\eta)) \in \ahu$ for all $x \in \ahu$.

\medskip

Choose $\th \in \ahu^*$ such that there exists $y \in \tilde{A}$ such that
$\th(\lambdau(\om)) = \tilde{\om}(y)$ for all $\om \in \duals$.

\medskip

Fix $\om \in \duals$ for the moment. Proposition \ref{gen.prop1} implies that
\begin{eqnarray*}
& & \th\bigl((\om_{v,w} \ot \io)(s_\mu(\lambdau(\om))\, \mu(\eta))\bigr)
= \th\bigl((\om_{v,w} \ot \io)(\mu(\om \eta))\bigr) \\
& & \spat = \th\bigl(\lambdau((\om \eta)\,\lambda^*(\om_{v,w}))\bigr)
= (\om \eta)\,\tilde{}\,\,(\lambda^*(\om_{v,w})\,y) \ .
\end{eqnarray*}
By the remarks before this proposition, we have for all $z \in A$ that
$$(\om \eta)\,\tilde{}\,\,(\lambda^*(\om_{v,w})\,z)
= (\om \eta)(\lambda^*(\om_{v,w})\,z)
= \langle \xi(\om \eta) , z^*  \la(\lambda^*(\om_{v,w})) \rangle\ .$$
Hence Kaplansky's density theorem implies that
$$\th\bigl((\om_{v,w} \ot \io)(s_\mu(\lambdau(\om))\, \mu(\eta))\bigr) 
= (\om \eta)\,\tilde{}\,\,(\lambda^*(\om_{v,w})\,y) = \langle \xi(\om \eta) , y^*  \la(\lambda^*(\om_{v,w})) \rangle\ .$$
So we get that 
$$\th\bigl((\om_{v,w} \ot \io)(s_\mu(\lambdau(\om))\, \mu(\eta))\bigr)
= \langle \lambda(\om) \, \xi(\eta) , y^* \la(\lambda^*(\om_{v,w}))\rangle
= \langle s_{\lambda}(\lambdau(\om))\,  \xi(\eta) , y^* \la(\lambda^*(\om_{v,w}))\rangle$$
Therefore 
$$\th\bigl((\om_{v,w} \ot \io)(s_\mu(p) \, \mu(\eta))\bigr) = 
\langle s_{\lambda}(p)\,  \xi(\eta) , y^* \la(\lambda^*(\om_{v,w}))\rangle \ .$$

Because $s_{\lambda}(p) = 0$, the previous equality implies that
$\th\bigl((\om_{v,w} \ot \io)(s_\mu(p)\,\mu(\eta))\bigr) = 0$. Therefore corollary \ref{gen.cor1} implies that $(\om_{v,w} \ot \io)(s_\mu(p)\,\mu(\eta))= 0$.
Since such elements $v$ and $w$  span a dense subset of $H$,  we conclude that $s_\mu(p) \, \mu(\eta) = 0$. So the density of  $\mu(\cI \cap L^1_*(A))$  in $s_\mu(\ahu)$ (see lemma 8.33 of \cite{J-V}) gives $s_\mu(p) = 0$.
\end{demo}

\medskip

Define  $K$ to be the closure of the subspace $\mu(\duals)(H \otimes H_{\text{\tiny u}})$ in $H \ot H_{\text{\tiny u}}$. In the following, $P$ will denote the  orthogonal projection on $K$.

\begin{corollary}
There exists a unique partial isometry $U \in M(A \ot B_0(H \ot H_{\text{\tiny u}}))$ such that $U^* U = U U^* = 1 \ot P$ and $\mu(\om) = (\om \ot \io)(U)$ for all $\om \in \duals$.
\end{corollary}
\begin{demo}
Since $\ker s_\lambda \subseteq \ker s_\mu$ and $s_\lambda(\ahu) = \ah$, there exists a unique $^*$-homomorphism $\phi : \ah \rightarrow B(H \ot H_{\text{\tiny u}})$ such that $\phi(s_\lambda(x)) = s_\mu(x)$ for all $x \in \ahu$. So $\phi(\lambda(\om)) = \mu(\om)$ for all $\om \in \duals$. It is also clear that $K$ is the closure of the subspace $\phi(\ah)(H \ot H_{\text{\tiny u}})$.

\medskip

By proposition 5.8 of \cite {Lan}, we know that $\phi$ is strict from $\ah$ to $B(H \ot H_{\text{\tiny u}})$ and $\phi(1) = P$. Remembering that $W \in M(A \ot \ah)$, we define
$U = (\io \ot \phi)(W)$. Then $U$ is a partial isometry in $M(A \ot B_0(H \ot H_{\text{\tiny u}}))$ such that
$U^* U = U U^* = (\io \ot \phi)(1) = 1 \ot P$. We have also for every $\om \in \duals$ that 
$$(\om \ot \io)(U) = (\om \ot \io)((\io \ot \phi)(W)) = \phi((\om \ot \io)(W)) = \phi(\lambda(\om)) = \mu(\om) \ .$$
\end{demo}

By the characterization of $U$ in the previous corollary and remark \ref{gen.rem2}, we have for all $\om \in \duals$ and $\eta \in B(H)_*$ that 
$$(\tilde{\om} \ovt \eta \ovt \io)(U) = (\eta \ovt \io)((\om \ot \io)(U))= (\eta \ovt \io)(\mu(\om)) \in \ahu\ .$$
Since $\duals$ is dense  in $\tilde{A}_*$, we conclude that
 $(\rho \ovt \io)(U) \in \ahu$ for all $\rho \in (\tilde{A} \ovt B(H))_*$.
Thanks to this simple observation and corollary \ref{gen.cor1}, we can copy the proof of theorem 3.1.4 in \cite{E}.

\begin {proposition} \label{gen.prop2}
There exists a unique element $\hat{\cV} \in M(A \ot B_0(H_{\text{\tiny u}}))$ such that
$\lambdau(\om) = (\om \ot \io)(\hat{\cV})$ for all $\om \in \duals$. 
Furthermore, $\hat{\cV}$ is a unitary element in $M(A \ot \ahu)$ such that 
$(\de \ot \io)(\hat{\cV}) = \hat{\cV}_{13} \hat{\cV}_{23}$.
\end{proposition}
\begin{demo}
Take $\th \in \ahu^*$ such that there exists $y \in \tilde{A}$ satisfying
$\th(\lambdau(\om)) = \tilde{\om}(y)$ for all $\om \in \duals$.

Then we have for all $\om \in \duals$ and $\eta \in B(H)_*$ that
\begin{eqnarray*}
\th\bigl((\tilde{\om} \ovt \eta \ovt \io)(U)\bigr)
& = & \th\bigl((\eta \ovt \io)(\mu(\om))\bigr) 
= \th\bigl(\lambdau(\om \,(\io \ovt \eta)(W))\bigr) \\
& = & \tilde{\om}((\io \ovt \eta)(W)\,y)
= (\tilde{\om} \ovt \eta)(W(y \ot 1)) 
\end{eqnarray*}
Since $\duals$ is dense in $\tilde{A}_*$, we get that $\th((\rho \ovt \io)(U)) = \rho(W(y \ot 1))$ for all $\rho \in B(H \ot H)_*$. Hence $\th((\rho \ovt \io)(W^*_{12}\,U)) = \rho(y \ot 1)$ for all $\rho \in 
B(H \ot H)_*$.

\smallskip

In particular, we have for all $\om \in \duals$ and $\eta \in B(H)_*$ that
\begin{equation} \label{gen.eq3}
\th\bigl((\tilde{\om} \ovt \eta \ovt \io)(W^*_{12}\, U)\bigr)
= \tilde{\om}(y)  \, \eta(1) = \th(\lambdau(\om)) \, \eta(1) \ .
\end{equation}

Therefore corollary \ref{gen.cor1} implies for every $\om \in \duals$ and $\eta \in B(H)_*$ that $(\tilde{\om} \ovt \eta \ovt \io)(W^*_{12}\,U) = \lambdau(\om) \, \eta(1)$. Hence 
$$(\om \ot \io \ot \io)(W^*_{12}\, U) = 1 \ot \lambdau(\om)$$
for all $\om \in \duals$.

\medskip

Equation \ref{gen.eq3} implies immediately that $W^*_{12}\, U (1 \ot x \ot 1) = (1 \ot x \ot 1)W^*_{12}\, U$ for all $x \in B(H)$. Therefore $W_{12}^*\, U \in (1 \ot B(H) \ot 1)'$ which implies the existence of $\hat{\cV} \in B(H \ot H_{\text{\tiny u}})$ such that $W_{12}^*\,U = \hat{\cV}_{13}$. Since $W \in M(A \ot B_0(H))$ and $U \in M(A \ot B_0(H) \ot B_0(H_{\text{\tiny u}}))$, we have clearly that $\hat{\cV} \in M(A \ot B_0(H_{\text{\tiny u}}))$. Of course, equation \ref{gen.eq3} also implies that $\lambdau(\om) = (\om \ot \io)(\hat{\cV})$ for all $\om \in \duals$.

\medskip

Moreover, the unitarity of $W$ implies that 
$$1 \ot P = U^* U = U^* W_{12}\, W_{12}^* \,U = \hat{\cV}_{13}^*\, \hat{\cV}_{13}\ .$$
This implies the existence of a projection $Q \in B(H_{\text{\tiny u}})$ such that $1 \ot Q = P$. Then we have also that $\hat{\cV} \hat{\cV}^* = 1 \ot Q$. Furthermore,
\begin{eqnarray*}
\hat{\cV}_{13} \,\hat{\cV}_{13}^* & = & W_{12}^* \,U U^* W_{12} = W_{12}^*\, (1 \ot P) W_{12} = W_{12}\, (1 \ot 1 \ot Q) W_{12}^* \\
& = & W_{12}\, W_{12}^*\, (1 \ot 1 \ot Q) = 1 \ot 1 \ot Q \ , 
\end{eqnarray*}
thus $\hat{\cV} \hat{\cV}^* = 1 \ot Q$. So $\hat{\cV}$ is a partial isometry with initial and final projection $1 \ot Q$.

\medskip

Take $v \in H_{\text{\tiny u}}$ such that $Q \,v = 0$. Then we have for every $w \in H$ that $(1 \ot Q)(w \ot v)  $ and hence  $\hat{\cV} (w \ot v) = 0$ (since $1 \ot Q$ is the initial projection of $\hat{\cV}$). This implies  for every $\om \in B_0(H)^*$ that $(\om \ot \io)(\hat{\cV})\, v = 0$. Consequently, $\lambdau(\om)\, v = 0$ for all $\om \in \duals$. Therefore the non-degeneracy of $\lambdau$ implies that $v = 0$. We conclude that $Q=1$ thus $\hat{\cV}$ is unitary.

\medskip

We have for all $\om_1,\om_2 \in \duals$ that
\begin{eqnarray*}
& & (\om_1 \ot \om_2 \ot \io)(\hat{\cV}_{13} \hat{\cV}_{23}) = (\om_1 \ot \io)(\hat{\cV}) \, (\om_2 \ot \io)(\hat{\cV}) = \lambda(\om_1) \, \lambda(\om_2) \\
& & \spat  = \lambda(\om_1 \om_2) = (\om_1 \om_2 \ot \io)(\hat{\cV})
= (\om_1 \ot \om_2 \ot \io)((\de \ot \io)(\hat{\cV})) \ .
\end{eqnarray*}

As usual, this implies that $(\de \ot \io)(\hat{\cV}) = \hat{\cV}_{13} \hat{\cV}_{23}$. Theorem 1.6 of \cite{Wor5} implies  that $\hat{\cV}$ is a unitary element in $M(A \ot \ahu)$.
\end{demo}

This proposition implies immediately that every $^*$-representation, and not only $\lambdau$, has a generator:

\begin{corollary} \label{gen.cor2}
Consider a \cst-algebra $C$ and a non-degenerate $^*$-homomorphism $\th : \duals \rightarrow M(C)$. Then there exists a unique element $V \in M(A \ot C)$ such that
$\th(\om) = (\om \ot \io)(V)$ for all $\om \in \duals$.
We have moreover that $V$ is unitary and $(\de \ot \io)(V) = V_{13} V_{23}$.
\end{corollary}
\begin{demo}
By the universal property of $\ahu$, there exists a unique non-degenerate $^*$-homomorphism 
$\bar{\th} : \ahu \rightarrow M(C)$ such that $\bar{\th} \lambdau = \th$. Now put 
$V = (\io \ot \bar{\th})(\hat{\cV})$.
\end{demo}

\begin{proposition} \label{gen.prop3}
Consider a \cst-algebra $C$ and a unitary element $V \in M(A \ot C)$ such that \newline $(\de \ot \io)(V) = V_{13} V_{23}$. Then there exists a unique non-degenerate $^*$-homomorphism $\th : \ahu \rightarrow M(C)$ such that $(\io \ot \th)(\hat{\cV}) = V$.
\end{proposition}
\begin{demo}
Uniqueness follows from proposition \ref{gen.prop2}. Define the linear mapping $\kappa : \duals \rightarrow M(C) : \om \mapsto (\om \ot \io)(V)$. A calculation like in the last part of the proof of proposition \ref{gen.prop2} shows that $\kappa$ is multiplicative.

\medskip

In the next part, we show that $\kappa$ is self adjoint. Choose $\eta \in C^*$. Let $\rho \in B_0(H)^*$ and define $\Upsilon \in B_0(H)^*$ such that $\Upsilon(x) = (\rho \ot \eta)(V (x \ot 1) V^*)$ for all $x \in B_0(H)$.

\smallskip

By assumption, $V_{13}\, V_{23} = (\de \ot \io)(V) = W_{12}^*\, V_{23}\, W_{12}  $, thus
$ W_{12} \, V_{13} = V_{23} \, W_{12} \, V_{23}^* $. \inlabel{gen.eq4}

Applying $\io \ot \rho \ot \eta$ to this equation gives $(\io \ot \rho)(W)\,(\io \ot \eta)(V) = (\io \ot \Upsilon)(W)$. Therefore proposition 8.3 of \cite{J-V} implies that $(\io \ot \rho)(W)\,(\io \ot \eta)(V) \in D(S)$ and 
\begin{eqnarray*}
& & S\bigl((\io \ot \rho)(W)\,(\io \ot \eta)(V)\bigr) = (\io \ot \Upsilon)(W^*)
= (\io \ot \rho \ot \eta)(V_{23} \, W_{12}^* \,V_{23}^*) \\
& & \spat \stackrel{(*)}{=} (\io \ot \rho \ot \eta)(V_{13}^* \, W_{12}^*)
= (\io \ot \eta)(V^*) \, (\io \ot \rho)(W^*) = (\io \ot \eta)(V^*)\, S((\io \ot \rho)(W)) \ ,
\end{eqnarray*}
where we used the adjoint of equation \ref{gen.eq4} in (*). Because such elements $(\io \ot \rho)(W)$ form a core for $S$, the closedness of $S$ implies for every $x \in D(S)$ that $x \, (\io \ot \eta)(V) \in D(S)$ and $S(x\,(\io \ot \eta)(V)) = (\io \ot \eta)(V^*)\,S(x)$.
From this, we infer that $(\io \ot \eta)(V) \in D(\bar{S})$ and 
$S((\io \ot \eta)(V)) = (\io \ot \eta)(V^*)$ \ (see remark 5.44 of \cite{J-V}).

\medskip

Choose $\om \in \duals$. By definition of $\om^*$, we have that $\om^*(x) = \overline{\om}(S(x))$ for all $x \in D(S)$. Since $D(S)$ is a strict \lq bounded\rq\ core for $\bar{S}$ (see remark 5.44 of \cite{J-V}), this gives $\om^*(x) = \overline{\om}(S(x))$ for all $x \in D(\bar{S})$. By the discussion above, we get for all $\eta \in C^*$ that
\begin{eqnarray*}
& & \eta((\om^* \ot \io)(V)) = \om^*((\io \ot \eta)(V)) = \overline{\om}\bigl(S((\io \ot \eta)(V))\bigr) \\
& & \spat = \overline{\om}((\io \ot \eta)(V^*)) = \eta((\overline{\om} \ot \io)(V^*))
= \eta((\om \ot \io)(V)^*) \ .
\end{eqnarray*} 
Consequently, $\kappa(\om^*) = \kappa(\om)^*$.

\medskip

Let us also verify quickly that $\kappa$ is non-degenerate: 
\begin{eqnarray*}
\overline{\kappa(\duals)\,C} & = & [\,(\om \ot \io)(V(1 \ot c)) \mid \om \in \duals, c \in C \,] \\
& = & [\,(\om \ot \io)(V(1 \ot c)) \mid \om \in \dual , c \in C \, ] \\
& \supseteq & [\,(a\,\om \ot \io)(V(1 \ot c)) \mid \om \in \dual , a \in A, c \in C \, ] \\
& = & [\, (\om \ot \io)(V(a \ot c)) \mid \om \in \dual , a \in A , c \in C \, ] \ .
\end{eqnarray*}
Because $V$ is unitary, this gives
$$ \overline{\kappa(\duals)\,C} \supseteq [\, (\om \ot \io)(b \ot d) \mid  b \in A , d \in C \, ] = C \ . $$

So we conclude that $\kappa : \duals \rightarrow M(C)$ is a non-degenerate $^*$-homomorphism. Therefore the universal property of $\ahu$ implies the existence of a non-degenerate $^*$-homomorphism $\th : \ahu \rightarrow M(C)$ such that $\th \lambdau = \lambda$.
We have for every $\om \in \duals$ that 
$$(\om \ot \io)((\io \ot \th)(\hat{\cV})) = \th((\om \ot \io)(\hat{\cV})) = \th(\lambdau(\om)) 
= \kappa(\om) = (\om \ot \io)(V) \ .$$
Hence $(\io \ot \th)(\hat{\cV}) = V$.
\end{demo}

\medskip

For later purposes we will need the projection from $\ah_{\text{\tiny u}}$ to $\ah$:

\begin{notation}
We define $\pih : \ahu \rightarrow \ah$   to be the surjective $^*$-homomorphism such that
$\pih \, \lambdau = \lambda$.
\end{notation}

We have for all $\om \in \duals$ that 
$$(\om \ot \io)\bigl((\io \ot \pih)(\hat{\cV})\bigr) =\pih\bigl((\om \ot \io)(\hat{\cV})\bigr)
= \pih(\lambdau(\om))  = \lambda(\om) = (\om \ot \io)(W) \ , $$
implying that $(\io \ot \pih)(\hat{\cV}) = W$.

\bigskip\bigskip

We can of course do the same thing for $L^1_*(\ah)$ and get the universal companion of $(A,\de)$ in this way.  So we define $\au$ to be the universal enveloping \cst-algebra of the Banach $^*$-algebra $L^1_*(\ah)$. 

By choosing the GNS-construction of the left Haar weight of $(\ah,\deh)$ in the right way, the multiplicative unitary of $(\ah,\deh)$ with respect to this well-chosen GNS-construction is equal to $\Sigma W^* \Sigma$ (see the remarks after proposition 8.20 of \cite{J-V}). This implies immediately that $(\ahh,\dhh) = (A,\de)$.
We will denote the  embedding of $L^1(\ah)$ into $A = \ahh$ by $\lambdah$. Notice that $\lambdah(\om) = (\io \ot \om)(W^*)$ for all $\om \in L^1(\ah)$. \inlabel{gen.eq8}

\smallskip

The embedding of $L^1(\ah)$ into $\au$ on the other hand will be denoted by $\lambdahu$. Define $\pi : \au \rightarrow A$ to be the surjective $^*$-homomorphism such that $\pi \, \lambdah_{\text{\tiny u}} = \lambdah$.

\medskip

As in proposition \ref{gen.prop2}, there exists a unitary element $\cV \in M(\au \ot \ah)$ such that  $\lambdau(\om) = (\io \ot \om)(\cV^*)$ for all $\om \in L^1_*(\ah)$.
We have moreover that $(\io \ot \deh)(\cV) = \cV_{13} \cV_{12}$ and $(\pi \ot \io)(\cV) = W$.
\inlabel{gen.eq5}

\smallskip

Also notice that
\begin{equation} \label{gen.eq6}
\au = [ \, (\io \ot \om)(\cV) \mid \om \in   L^1(\ah)\, ] = [ \, (\io \ot \om)(\cV) \mid \om \in   B_0(H)^*\, ]  \ .
\end{equation}

\medskip

In this setting, proposition \ref{gen.prop3} gets the following form.

\begin{proposition} \label{gen.prop4}
Consider a \cst-algebra $C$ and a unitary element $V \in M(C \ot \ah)$ such that \newline $(\io \ot \deh)(V) = V_{13} V_{12}$. Then there exists a unique non-degenerate $^*$-homomorphism $\th : \au \rightarrow M(C)$ such that $(\th \ot \io)(\cV) = V$.
\end{proposition}

\medskip

Although $\au$ is defined to be the universal enveloping \cst-algebra of a space of linear functionals on $\ah$, it is also the universal enveloping \cst-algebra of a dense subalgebra of $A$. Consider the  injective algebra homomorphism $\lambdah :  L^1(\ah) \rightarrow A $.

Proposition 8.32 of \cite{J-V} implies that $\lambda(L^1_*(\ah)) = \lambda(L^1(\ah)) \cap \lambda(L^1(\ah))^*$ and that the restriction of $\lambda$ to $L^1_*(\ah)$ is a $^*$-isomorphism from $L^1_*(\ah)$ to $\lambda(L^1(\ah)) \cap \lambda(L^1(\ah))^*$.

This implies the following. Define $\cA = \{\,(\io \ot \om)(W)\mid \om \in B_0(H)^*\,\}$. Then $\cA$ is a dense subalgebra of $A$, $\cA \cap \cA^*$ is a dense sub $^*$-algebra of $A$ and 
$\au$ is the universal enveloping \cst-algebra of $\cA \cap \cA^*$.

\medskip

A similar remark applies to $\ah$ but in this case we have to replace the algebra $\cA$ by the algebra $\hat{\cA}$ defined by $\hat{\cA} = \{\,(\om \ot \io)(W)\mid \om \in B_0(H)^*\,\}$.

\bigskip\medskip

\sectie{The universal bi-\cst-algebras,  the universal corepresentation}

Up to now, we only have constructed \lq universal\rq\  \cst-algebras $\au$ and $\ahu$. In this section, we introduce the comultiplications on them and construct the universal corepresentation between them.

\medskip

First we repeat a standard terminology in quantum group theory.

\begin{terminology}
Consider a bi-\cst-algebra $(B,\de)$, a \cst-algebra $C$ and a unitary element $V \in M(B \ot C)$ such that $(\de \ot \io)(V) = V_{13} V_{23}$. Then $V$ is called a unitary corepresentation of $(B,\de)$ on $C$.
\end{terminology} 

\medskip

In a first step, we will follow the standard road to introduce the comultiplication and counit on $\au$ (see theorem 1.3 of \cite{PW}). We want to define the comultiplication $\deu$ on $\au$ in such a way that $\cV$ is a unitary corepresentation of $(\au,\deu)$.

\begin{proposition} \label{com.prop4}
There exists a unique non-degenerate $^*$-homomorphisms $\deu : \au \rightarrow M(\au \ot \au)$
such that $(\deu \ot \io)(\cV) = \cV_{13} \cV_{23}$. We have moreover that
\begin{enumerate}
\item $(\deu \ot \io)\deu = (\io \ot \deu)\deu$.
\item $\deu(\au)(\au \ot 1)$ and $\deu(\au)(1 \ot \au)$ are dense subspaces in $\au \ot \au$.
\end{enumerate}
\end{proposition}
\begin{demo}
We want to apply proposition \ref{gen.prop4} in order to get hold of $\deu$. Therefore look at the unitary element $\cV_{13} \cV_{23} \in M(\au \ot \au \ot \ah)$. We have that \begin{eqnarray*}
& & (\io \ot \io \ot \deh)(\cV_{13} \cV_{23}) = (\io \ot \io \ot \deh)(\cV_{13})
\, (\io \ot \io \ot \deh)(\cV_{23}) \\
& & \spat \stackrel{(*)}{=} \cV_{14}  \cV_{13} \cV_{24} \cV_{23}
= \cV_{14} \cV_{24} \cV_{13} \cV_{23}
= \bigl(\cV_{13} \cV_{23}\bigr)_{13} \, \bigl(\cV_{13} \cV_{23}\bigr)_{12} \ ,
\end{eqnarray*}
where we used the first equality of equation \ref{gen.eq5} in equality (*).
Therefore proposition \ref{gen.prop4} implies the existence of a unique non-degenerate $^*$-homomorphism $\deu$ from 
$\au$ to $M(\au \ot \au)$ such that $(\deu \ot \io)(\cV) = \cV_{13} \cV_{23}$.

\medskip

We have that
$$((\deu \ot \io)\deu \ot \io)(\cV)
= (\deu \ot \io \ot \io)(\cV_{13} \cV_{23}) = \cV_{14} \cV_{24} \cV_{34}$$
and 
$$((\io  \ot \deu)\deu \ot \io)(\cV)
= ( \io \ot \deu \ot \io)(\cV_{13} \cV_{23})
= \cV_{14} \cV_{24} \cV_{34} \ , $$
implying that
$$((\deu \ot \io)\deu \ot \io)(\cV) =  ((\io  \ot \deu)\deu \ot \io)(\cV) \ .$$
Therefore $(\deu \ot \io)\deu = (\io  \ot \deu)\deu$ by equation \ref{gen.eq6}.

\medskip

Let us now verify the density conditions (see proposition 5.1 of \cite{Wor5}). Using equation \ref{gen.eq6}, we see that
\begin{eqnarray*}
\overline{\deu(\au)(1 \ot \au)}
& = & [ \, \deu((\io \ot \om)(\cV))(\,1 \ot a)\mid \om \in B_0(H)^*,a \in \au \,] \\
& = & [ \, (\io \ot \io \ot \om)(\cV_{13} \cV_{23} (1 \ot a \ot 1)) 
\mid \om \in B_0(H)^* ,a \in \au \,] \\
& = & [ \, (\io \ot \io \ot x \om)(\cV_{13} \cV_{23} (1 \ot a \ot 1)) 
\mid \om \in B_0(H)^*, x \in B_0(H),  a \in \au \,]  \\
& = & [ \, ( \io \ot \io \ot \om)(\cV_{13} \cV_{23} (1 \ot a \ot x)) 
\mid \om \in B_0(H)^*, x \in B_0(H),  a \in \au \,]\ .
\end{eqnarray*}
Since $\cV$ is a unitary element in $M(\au \ot B_0(H))$, we have that 
$\cV(\au \ot B_0(H)) = \au \ot B_0(H)$.
Hence the above chain of equalities implies that 
\begin{eqnarray*}
\overline{\deu(\au)(1 \ot \au)}
& = & [ \, (\io \ot \io \ot \om)(\cV_{13} (1 \ot a \ot x)) 
\mid \om \in B_0(H)^*, x \in B_0(H),  a \in \au \,] \\
& = & [ \, (\io \ot \io \ot x \, \om)(\cV_{13})\,(1 \ot a) \mid \om \in B_0(H)^*, x \in B_0(H),  a \in \au \,] \\
& = & [\, (\io \ot \io \ot \om)(\cV_{13})\,(1 \ot a) \mid \om \in B_0(H)^*,   a \in \au \,]\\
& = & [\, (\io \ot \om)(\cV) \ot a \mid \om \in B_0(H)^*,   a \in \au \,]
= \au \ot \au \ .
\end{eqnarray*}
In a similar way, one proves that $(\au \ot 1)\deu(\au)$ is a dense subspace of $\au \ot \au$,
and we are done (remember the $^*$-operation).
\end{demo}

\begin{proposition} \label{com.prop1}
The following identities hold:
\begin{enumerate}
\item $(\pi \ot \pi)\deu = \de \pi$,
\item $(\io \ot \pi)(\deu(x)) = \cV^*\, (1 \ot \pi(x))\, \cV$ for all $x \in \au$.
\end{enumerate}
\end{proposition}
\begin{demo}
Since $\de$ is implemented by $W$ and $(\pi \ot \io)(\cV) = W$, the first equality follows from the second one. So we only have to prove the second statement.

\smallskip

Using the second equality in equation \ref{gen.eq5}, we get that 
$$((\io \ot \pi)\deu \ot \io)(\cV)
= (\io \ot \pi \ot \io)(\cV_{13} \cV_{23})  = \cV_{13} W_{23}$$
whereas
$$\cV_{12}^* (1 \ot (\pi \ot \io)(\cV)) \cV_{12}
= \cV_{12}^* W_{23} \cV_{12}$$
Since $\cV_{12} \cV_{13} = (\io \ot \flip \deh)(\cV) = W_{23} \cV_{12} W_{23}^*$, we conclude that
$$((\io \ot \pi)\deu \ot \io)(\cV) = \cV_{12}^* (1 \ot (\pi \ot \io)(\cV)) \cV_{12}$$
So we get for all $\om \in B_0(H)^*$ that 
$$(\io \ot \pi)\bigl(\deu((\io \ot \om)(\cV))\bigr)
= \cV^* (1 \ot \pi((\io \ot \om)(\cV))\,) \cV  $$
and statement 2. follows from equation \ref{gen.eq6}.
\end{demo}

\medskip

If we apply proposition \ref{gen.prop4} to the unit in $M(\au \ot \ah)$, we get hold of the counit $\vepu$ on $(\au,\deu)$.

\begin{proposition} \label{com.prop3}
There exists a unique non-zero *-homomorphism $\vepu : \au \rightarrow \C$  such that 
$$(\vepu \ot \io)\deu = (\io \ot \vepu)\deu = \io \ .$$
Moreover, $(\vepu \ot \io)(\cV) = 1$.
\end{proposition}
\begin{demo}
Uniqueness is trivial. By proposition \ref{gen.prop4}, there exists a unique non-zero *-homomorphism $\vepu : \au \rightarrow \C$ such that
$(\vepu \ot \io)(\cV) = 1$. Therefore, 
$$((\vepu \ot \io)\deu \ot \io)(\cV) = (\vepu \ot \io \ot \io)(\cV_{13} \cV_{23})
= (1 \ot (\vepu \ot \io)(\cV)) \, \cV = \cV \ .$$
Hence equation \ref{gen.eq6} implies that $(\vepu \ot \io)\deu = \io$. Similarly, $(\io \ot \vepu)\deu = \io$.
\end{demo}

Notice that this proposition implies that $\deu : \au \rightarrow M(\au \ot \au)$ is injective.

\bigskip\medskip

Of course, we also have corresponding results for $\ahu$. Let us explicitly formulate them.

\begin{proposition}
There exists a unique non-degenerate $^*$-homomorphisms $\dehu : \ahu \rightarrow M(\ahu \ot \ahu)$
such that $(\io \ot \dehu)(\hat{\cV}) = \hat{\cV}_{13} \hat{\cV}_{12}$. We have moreover that
\begin{enumerate}
\item $(\dehu \ot \io)\dehu = (\io \ot \dehu)\dehu$.
\item $\dehu(\ahu)(\ahu \ot 1)$ and $\dehu(\ahu)(1 \ot \ahu)$ are dense subspaces of $\ahu \ot \ahu$.
\end{enumerate}
\end{proposition}

\medskip

\begin{proposition} \label{com.prop2}
The following identities hold:
\begin{enumerate}
\item We have that $(\pih \ot \pih)\dehu = \dehu \pih$.
\item $(\pih \ot \io)(\flip \dehu(x)) = \hat{\cV}\,(\pih(x) \ot 1)\, \hat{\cV}^*$
for all $x \in \ahu$.
\end{enumerate}
\end{proposition}

\bigskip\medskip

It turns out to be very easy to get hold of the \lq universal\rq\ corepresentation between $\au$ and $\ahu$. First we need a lemma which we copied from proposition 3.11 of \cite{Mas-Nak}.

\begin{lemma}
The set $M(A) \cap M(\ah)$ is equal to $\C\,1$.
\end{lemma}
\begin{demo}
Choose $x \in M(A) \cap M(\ah)$. By proposition 8.17 of \cite{J-V}, we know that $\hat{R}(x) = J x^* J$. Because $x$ belongs to $M(A)$, Tomita-Takesaki theory tells us that $J x^* J \in A'$, so $\hat{R}(x) \in A'$.
Thus, since $W \in M(A \ot \ah)$, 
$$\deh(\hat{R}(x)) = W (\hat{R}(x) \ot 1) W^* = (\hat{R}(x) \ot 1) W W^* = \hat{R}(x) \ot 1 \ .$$
Therefore result 6.1 of \cite{J-V} implies that $\hat{R}(x) \in \C \, 1$, so $ x \in \C\,1$.
\end{demo}

\medskip

The next result guarantees that the multiplicative unitary $W$ is basic in the sense of definition 2.3 of \cite{Ng} and that $\cU$ is universal in the sense of this same definition. See also lemma 1.5 of the same paper.

\begin{proposition}
There exists a unique unitary element $\cU \in M(\au \ot \ahu)$ such that
$\cU_{13} = \cV_{12}^* \hat{\cV}_{23} \cV_{12} \hat{\cV}_{23}^*$.
We have moreover that $(\deu \ot \io)(\cU) = \cU_{13}\, \cU_{23}$ and
$(\io \ot \dehu)(\cU) = \cU_{13} \, \cU_{12}$. The element  $\cU$ is called the universal corepresentation of $(\au,\deu)$.
\end{proposition}
\begin{demo}
Since $\hat{\cV} \in M(A \ot \ahu)$, proposition \ref{com.prop1}.2 implies that $\cV_{12}^* \hat{\cV}_{23} \cV_{12}$ belongs to $M(\au \ot A \ot \ahu)$. Hence $\cV_{12}^* \hat{\cV}_{23} \cV_{12} \hat{\cV}_{23}^*$ belongs to $M(\au \ot A \ot \ahu)$.
Because $\cV \in M(\au \ot \ah)$, proposition \ref{com.prop2}.2 implies that $\hat{\cV}_{23} \cV_{12} \hat{\cV}_{23}^*$ belongs to
$M(\au \ot \ah \ot \ahu)$. Hence $\cV_{12}^* \hat{\cV}_{23} \cV_{12} \hat{\cV}_{23}^*$ belongs to $M(\au \ot \ah \ot \ahu)$.

Therefore the previous lemma  implies the existence of a unitary element $\cU \in M(\au \ot \ahu)$ such that $\cU_{13} = \cV_{12}^* \hat{\cV}_{23} \cV_{12} \hat{\cV}_{23}^*$.

\medskip

We have that $\cV_{12}\,\cU_{13} = \hat{\cV}_{23} \cV_{12} \hat{\cV}_{23}^*$. Since $\cV$ is a corepresentation of $(\au,\deu)$ on $B_0(H)$, this equality implies easily that $\cV_{12} \, \cU_{13}$ is a corepresentation of $(\au,\deu)$ on $B_0(H) \ot \ah_u$.

Hence 
\begin{eqnarray*}
(\deu \ot \io \ot \io)(\cV_{12} \, \cU_{13}) & = & (\cV_{12} \, \cU_{13})_{13} \,(\cV_{12} \, \cU_{13})_{23}
= \cV_{13} \, \cU_{14} \cV_{23} \, \cU_{24}\\
& = & \cV_{13} \cV_{23} \,\cU_{14} \,\cU_{24}
= (\deu \ot \io \ot \io)(\cV_{12}) \, \cU_{14} \,\cU_{24} \ .
\end{eqnarray*}
Consequently, $(\deu \ot \io \ot \io)(\cU_{13}) = \cU_{14} \,\cU_{24}$.
Therefore $(\deu \ot \io)(\cU) = \cU_{13} \, \cU_{23}$. Similarly, $(\io \ot \dehu)(\cU) = \cU_{13} \, \cU_{12}$.
\end{demo}

\begin{corollary} \label{com.cor1}
We have that
\begin{enumerate}
\item $(\io \ot \pih)(\cU) = \cV$.
\item $(\pi \ot \io)(\cU) = \hat{\cV}$.
\item $(\pi \ot \pih)(\cU) = W$.
\end{enumerate}
\end{corollary}
\begin{demo}
\begin{enumerate}
\item Using the facts that $(\io \ot \pih)(\hat{\cV}) = W$ and $(\io \ot \flip \deh)(\cV) = W_{23} \cV_{12} W_{23}^*$, we get that
\begin{eqnarray*} 
(\io \ot \io \ot \pih)(\cU_{13}) & = & (\io \ot \io \ot \pih)(\cV_{12}^* \hat{\cV}_{23} \cV_{12} \hat{\cV}_{23}^*) 
= \cV_{12}^* W_{23} \cV_{12} W_{23}^* \\
& = & \cV_{12}^* (\io \ot \flip \deh)(\cV) = \cV_{12}^*  \cV_{12} \cV_{13} = \cV_{13} \ .
\end{eqnarray*}
Hence $(\io \ot \pih)(\cU) = \cV$.
\item Similar to the first equality.
\item Follows from the first result and the fact that $(\pi \ot \io)(\cV) = W$.
\end{enumerate}
\end{demo}

\begin{remark} \rm \label{com.rem1}
Notice that the previous result, equation \ref{gen.eq6} and its obvious dual version imply that
$$A_u \subseteq [\,(\io \ot \om)(\cU) \mid \om \in \ahu^* \, ] \hspace{1.5cm} \text{and} \hspace{1.5cm}
\ahu \subseteq [\,(\om \ot \io)(\cU) \mid \om \in \au^* \, ] \ .$$
\end{remark}

\bigskip

If $V$ is a unitary corepresentation of $(\au,\deu)$ on a \cst-algebra $C$, then proposition \ref{com.prop1}.1 implies that \newline $(\pi \ot \io)(V)$ is a unitary corepresentation of $(A,\de)$ on $C$. But it turns out that every unitary corepresentation of $(A,\de)$ lifts to a unique unitary corepresentation of $(\au,\deu)$ in this way.

\begin{result} \label{com.res1}
Consider a \cst-algebra $C$ and unitary corepresentations $U,V$ of $(\au,\deu)$ on $C$ such that $(\pi \ot \io)(U) = (\pi \ot \io)(V)$. Then $U = V$.
\end{result}
\begin{demo}
We have by assumption that $(\deu \ot \io)(U) = U_{13} U_{23}$. If we apply $\io \ot \pi \ot \io$ to this equality and use proposition \ref{com.prop1}.2, we get that $\cV_{12}^*\, (\pi \ot \io)(U)_{23}\, \cV_{12} = U_{13}\, (\pi \ot \io)(U)_{23}$ and therefore
$$U_{13} = \cV_{12}^* \, (\pi \ot \io)(U)_{23} \, \cV_{12} \, (\pi \ot \io)(U)_{23}^*\ .$$
Similarly, $V_{13} = \cV_{12}^* \, (\pi \ot \io)(V)_{23} \, \cV_{12} \, (\pi \ot \io)(V)_{23}^*$, so we get that $U_{13} = V_{13}$. Hence $U = V$.
\end{demo}

\medskip

Let us quickly explain the universal property of $\hat{\cU}$. The next proposition guarantees that $\cU$ induces a bijection between non-degenerate $^*$-homomorphisms of $\ahu$ and unitary corpresentations of $(\au,\deu)$. 

\begin{proposition}
Consider a \cst-algebra $C$ and a unitary corepresentation $U$ of $(\au,\deu)$ on $C$.  Then there exists a unique non-degenerate $^*$-homomorphism $\th : \ahu \rightarrow M(C)$ such that $(\io \ot \th)(\cU) = U$.
\end{proposition}
\begin{demo}
Uniqueness follows from remark \ref{com.rem1}. By proposition \ref{gen.prop3}, there exists a unique non-degenerate $^*$-homomorphism $\th : \ahu \rightarrow M(C)$ such that $(\io \ot \th)(\hat{\cV}) = (\pi \ot \io)(U)$. Therefore corollary \ref{com.cor1}.2 implies that
$$(\pi \ot \io)((\io \ot \th)(\cU)) = (\io \ot \th)(\hat{\cV}) = (\pi \ot \io)(U)  \ .$$
Hence, by the previous result, $(\io \ot \th)(\cU) = U$.
\end{demo}

\begin{proposition} \label{com.prop5}
Consider a \cst-algebra $C$ and a unitary corepresentation $U$ of $(A,\de)$ on $C$.  Then there exists a unique corepresentation $V$ of $(\au,\deu)$ on $C$ such that $(\pi \ot \io)(V) = U$.
\end{proposition}
\begin{demo}
Uniqueness follows from result \ref{com.res1}. By proposition \ref{gen.prop3} we get the existence of a  non-degenerate $^*$-homomorphism $\th : \ahu \rightarrow M(C)$ such that $(\io \ot \th)(\hat{\cV}) = U$. Put $V = (\io \ot \th)(\cU)$ which is a unitary corepresentation of $(\au,\deu)$ on $C$ such that 
$$(\pi \ot \io)(V) = (\pi \ot \io)((\io \ot \th)(\cU)) = (\io \ot \th)(\hat{\cV}) = U \ .$$
\end{demo}

So we have proven (in a very elementary way) that $(A,\de)$ and $(\au,\deu)$ have the same corepresentation theory.

\bigskip\medskip

\sectie{Lifting bi-automorphisms}

In the last part of the previous section, we showed that $(\au,\deu)$ and $(A,\de)$ have the same unitary corepresentation theory. The same is true for their bi-automorphisms, i.e. 
automorphisms which commute with the comultiplication. We will  work in a setting which is a little bit more general than the framework of bi-automorphisms in order to lift the modular groups of the Haar weights. In the last proposition of this section, we  lift the unitary antipode from the reduced to the universal level.

\bigskip

In the last statement of the next proposition, we use the language of Hilbert \cst-modules and the identification $M(B \ot B_0(H)) = \cL(B \ot H)$ for every \cst-algebra $B$.

\begin{proposition} \label{bi.prop1}
Consider  $^*$-automorphisms $\al$ and $\be$ on $A$ such that $(\al \ot \be)\de = \de \al$. Then the following properties hold.
\begin{enumerate}
\item $(\be \ot \be)\de = \de \be$.
\item There exists a number $r > 0$ such that $\vfi \, \al = r \, \vfi$ and $\vfi \, \be = r \, \vfi$.
\item Define two unitary operators $U$,$V$ on $H$ such that $U \la(a) = r^{-\frac{1}{2}}\, \la(\al(a))$ and $V \la(a) = r^{-\frac{1}{2}}\, \la(\be(a))$ for all $a \in \Nfi$. Then
$$(\al \ot \io)(W) = (1 \ot U^*) W (1 \ot V) \hspace{1cm} \text{ and } \hspace{1cm} (\be \ot \io)(W)  = (1 \ot V^*) W ( 1 \ot V) \ .$$
\item There exists unique $^*$-automorphisms $\alu$,$\beu$ on $\au$ such that
$$(\alu \ot \io)(\cV) = (1 \ot U^*) \cV (1 \ot V) \hspace{1cm} \text{ and } \hspace{1cm}
(\beu  \ot \io)(\cV)  = (1 \ot V^*) \cV ( 1 \ot V) \ .$$
We have moreover that $\pi \alu = \pi \al$, $\pi \beu = \be \pi$ and 
$(\alu \ot \beu)\deu = \deu \, \alu$.
\end{enumerate}
\end{proposition}
\begin{demo}
\begin{trivlist}
\item[\ \,1.] We have that
\begin{eqnarray*}
(\al \ot (\be \ot \be)\de)\de & = & (\al \ot \be \ot \be)(\de \ot \io)\de  
= (\de \ot \io)(\al \ot \be)\de \\
&  = & (\de \ot \io)\de \al = (\io \ot \de)\de \al = (\al \ot \de \be)\de \ .
\end{eqnarray*}
Hence $(\io \ot (\be \ot \be)\de)\de = (\io \ot \de \be)\de$. Therefore the density conditions in definition \ref{red.def1} imply that $(\be \ot \be)\de = \de \be$.
\item[\ \,2.] Since $(\be \ot \be)\de = \de \be$, the proper weight $\vfi  \be$ is  left invariant. Hence the uniqueness of the left Haar weight (see theorem 7.14 of \cite{J-V}) implies the existence of a number $r > 0$ such that $\vfi \be = r \, \vfi$.

Choose $a \in \Mfi^+$. Then we have by left invariance of $\vfi$ for all $\om \in A^*$ that
$(\om \al \ot \io)\de(a) \in \Mfi^+$ and therefore the relative invariance of $\vfi$ under $\be$  implies that $(\om \ot \io)\de(\al(a)) = \be((\om \al \ot \io)\de(a)) \in \Mfi^+$.
Consequently, proposition 6.2 of \cite{J-V} gives  that $\al(a) \in \Mfi^+$. Now take $\eta \in A^*_+$ such that $\eta(1) = 1$. Then the left invariance of $\vfi$ implies that
$$\vfi(\al(a)) = \vfi\bigl((\eta  \ot \io)\de(\al(a))\bigr)
= \vfi\bigl(\be((\eta \al \ot \io)\de(a))\bigr)
= r \, \vfi((\eta \al \ot \io)\de(a)) = r \, (\eta \al)(1) \, \vfi(a) = r \, \vfi(a)\ .$$

Working with $\al^{-1}$ instead of $\al$, also $\al^{-1}(\Mfi^+) \subseteq \Mfi^+$. Hence $\vfi \al = r \, \vfi$.

\item[\ \,3.] Take $\om \in L^1(A)$. Using result 2.10 of \cite{J-V}, we have for all $a \in \Nfi$ that 
\begin{eqnarray*}
& & (\om \al \ot \io)(W^*) \la(a) = \la((\om \al \ot \io)\de(a)) = \la\bigl(\be^{-1}((\om \ot \io)\de(\al(a)))\bigr) \\ 
& & \spat = r^{-\frac{1}{2}}\, V^* \la((\om \ot \io)\de(\al(a))) = r^{-\frac{1}{2}}\, V^* (\om \ot \io)(W^*) \, \la(\al(a)) = V^* (\om  \ot \io)(W^*) U \la(a) \ ,
\end{eqnarray*}
implying that $(\om \al \ot \io)(W^*) = V^* (\om \ot \io)(W^*) U$. So we conclude that  $(\al \ot \io)(W^*) = (1 \ot V^*) W^* (1 \ot U)$. 

Similarly, the relation $(\be \ot \be)\de = \de \be$ implies that $(\be \ot \io)(W)  = (1 \ot V^*) W ( 1 \ot V)$.

\item[\ \,4.] Uniqueness follows immediately from equation \ref{gen.eq6}, let us turn to the existence.

We could use proposition \ref{bi.prop1} to prove the existence but this is actually a detour. We  will illustrate the use of proposition \ref{bi.prop1} in the proof of the next result.
By the equalities in statement 3., we have for every $\om \in L^1(A)$ that 
$V^* (\om \ot \io)(W^*) U = (\om \al \ot \io)(W^*)$. Hence the definition of $\ah$ implies that $V^* \ah\, U = \ah$. Since $\cV \in M(\au \ot \ah)$, this implies that $(1 \ot V^*)\cV^*(1 \ot U)$ belongs to $M(A \ot \ah)$.

\smallskip

Define the linear  map $\tilde{\al} :  L_1(\ah) \rightarrow L_1(\ah)$ such that $\tilde{\al}(\om)(x) = \om(V^* x U)$ for all $\om \in L_1(\ah)$ and $x \in \ah$. Formula \ref{gen.eq8} implies that $\al(\lambdah(\om)) = \lambdah(\tilde{\al}(\om))$ for all $\om \in L^1(\ah)$. Because $\al$ is multiplicative, it follows easily that $\tilde{\al}$ is multiplicative. Since $\al$ is a self adjoint mapping, proposition 8.32 of \cite{J-V} implies that $\tilde{\al}(L^1_*(\ah)) \subseteq   L^1_*(\ah)$ and that $\tilde{\al}(\om)^* = \tilde{\al}(\om^*)$ for $\om \in L^1_*(\ah)$.

\smallskip

So the restriction of $\tilde{\al}$ to $L^1_*(\ah)$ is a $^*$-homomorphism from $L^1_*(\ah)$ into $L^1_*(\ah)$. Therefore the  universal property of $\au$ implies the existence of a $^*$-homomorphism $\alu : \au \rightarrow \au$ such that $\alu(\lambdau(\om)) = \lambdau(\tilde{\al}(\om))$ for all $\om \in L_1^*(\ah)$.
This implies for every $\om \in L^1_*(\ah)$ 
\begin{eqnarray*}
(\io \ot \om)((\alu \ot \io)(\cV^*)) & = & \alu((\io \ot \om)(\cV^*)) 
= \alu(\lambdau(\om)) = \lambdau(\tilde{\al}(\om)) \\
& = & (\io \ot \tilde{\al}(\om))(\cV^*) = (\io \ot \om)((1 \ot V^*)\cV^*(1 \ot U)) \ .
\end{eqnarray*}
Hence $(\alu \ot \io)(\cV^*) = (1 \ot V^*)\cV^*(1 \ot U)$.



\smallskip

We can of course do the same thing for $\al^{-1}$. This gives  a $^*$-homomorphism $\gammau : \au \rightarrow \au$ such that $(\gammau \ot \io)(\cV) = (1 \ot U) \cV (1 \ot V^*)$.
Then it is clear that $(\gammau \alu \ot \io)(\cV) = (\alu \gammau \ot \io)(\cV) = \cV$. From equation \ref{gen.eq6}, we conclude that  $\gammau \alu = \alu \gammau = \io$. Thus $\alu$ is an $^*$-automorphism on $\au$. Moreover, 
$$ (\pi \alu \ot \io)(\cV)  =  (\pi \ot \io)((1 \ot U^*)\cV(1 \ot V))
= (1 \ot U^*) W (1 \ot V) = (\al \ot \io)(W) = (\al \pi \ot \io)(\cV) \ ,$$
hence $\alu \pi = \pi \al$ by equation \ref{gen.eq6}. 

\smallskip

The $^*$-automorphism $\beu$ is constructed in a similar way. Moreover,  
\begin{eqnarray*}
& & ((\alu \ot \beu)\deu \ot \io)(\cV) 
= (\alu \ot \beu \ot \io)(\cV_{13} \cV_{23}) \\
& & \spat = [(1 \ot 1 \ot U^*) \cV_{13} (1 \ot 1 \ot V)]\,[(1 \ot 1 \ot V^*) \cV_{23} (1 \ot 1 \ot V)]
= (1 \ot 1 \ot U^*) \cV_{13} \cV_{23} (1 \ot 1 \ot V) \\
& & \spat = (1 \ot 1 \ot U^*)(\deu \ot \io)(\cV) (1 \ot 1 \ot V)
= (\deu \ot \io)((1 \ot U^*)\cV(1 \ot V)) = (\deu \alu \ot \io)(\cV) 
\end{eqnarray*}
and the last equation of the proposition  follows.
\end{trivlist}
\end{demo}

\bigskip

We want to use the same principle to lift the unitary antipode to the level of $\au$.

\begin{proposition} \label{bi.prop2}
There exists a unique $^*$-antiautomorphism $\Ru$ on $\au$ such that
$(\Ru \ot \hat{R})(\cV) = \cV$. We have moreover that $\Ru^2 = \io$, $\flip(\Ru \ot \Ru)\deu = \deu \Ru$ and  $\pi \Ru = R \pi$.
\end{proposition}
\begin{demo}
Denote the opposite \cst-algebra of $\au$ by $\au^\circ$ and let $\th : \au \rightarrow \au^\circ$ be the obvious $^*$-antiisomorphism. Then $(\th \ot \hat{R})(\cV)$ is a unitary element in $M(\au^\circ \ot \ah)$ such that
\begin{eqnarray*}
(\io \ot \deh)\bigl((\th \ot \hat{R})(\cV)\bigr) & = &  
(\io \ot \flip)(\th \ot \hat{R} \ot \hat{R})\bigl((\io \ot \deh)(\cV)\bigr)  
 = (\io \ot \flip)(\th \ot \hat{R} \ot \hat{R})(\cV_{13} \cV_{12})\\
& = & (\io \ot \flip)((\th \ot \hat{R})(\cV)_{12} (\th \ot \hat{R})(\cV)_{13})
= (\th \ot \hat{R})(\cV)_{13} (\th \ot \hat{R})(\cV)_{12} \ .
\end{eqnarray*}
Therefore proposition \ref{gen.prop4} guarantees the existence of a non-degenerate $^*$-homomorphism $\eta : \au \mapsto M(\au^\circ)$ such that $(\eta \ot \io)(\cV) = (\th \ot \hat{R})(\cV)$. Define $\Ru = \th^{-1} \, \eta$, so $\Ru$ is a non-degenerate $^*$-antihomomor\-phism from $\au$ to $M(\au)$
such that  $(\Ru \ot \hat{R})(\cV) = \cV$. By equation \ref{gen.eq6}, this equation implies that $\Ru(\au) = \au$.
Because $\hat{R}^2 = \io$, we also conclude that $(\Ru^2 \ot \io)\cV = \cV$ which implies that $\Ru^2 = \io$. We have also that
$$(\pi \Ru \ot \hat{R})\cV = (\pi \ot \io)\cV = W = (R \ot \hat{R})(W) = (R \pi \ot \hat{R})(\cV) \ ,$$
implying that $\pi \Ru = R \pi$.
\end{demo}

\bigskip\medskip

\sectie{Left and right Haar weights of the universal quantum group}

We use the surjective $^*$-homomorphism $\pi : \au \rightarrow A$ to pull back the left and right Haar weights on $(A,\de)$ to left and right invariant weights on $(\au,\deu)$. We prove that $\cV$ is the multiplicative unitary naturally associated to the resulting left invariant weight on $(\au,\deu)$. In the last part, a converse of proposition \ref{bi.prop1} is formulated.

\begin{proposition}
We define $\psiu = \psi \, \pi$. Then $\psiu$ is a proper weight on $\au$ which has a GNS-construction $(H,\pi,\gau)$ such that $\gau = \ga \, \pi$.
\end{proposition}
\begin{demo}
Choose $a,b \in \Nfi$ and $c \in \Npsi$. 
Then 
$$\pi\bigl((\io \ot \om_{\la(a), \la(c^* b)})(\cV)\bigr)
= (\io \ot \om_{\la(a),\la(c^* b)})(W)
= (\io \ot \vfi)(\de(b^* c)(1 \ot a))$$
By result 2.6 of \cite{J-V}, we know that $(\io \ot \vfi)(\de(b^* c)(1 \ot a))$ belongs to $\Npsi$, hence $\pi\bigl((\io \ot \om_{\la(a), \la(c b)})(\cV)\bigr)$ belongs to $\Nfi$.
We conclude that $(\io \ot \om_{\la(a), \la(c b)})(\cV)$ belongs to $\cN_{\vfiu}$.

\smallskip

Therefore equation \ref{gen.eq6} implies that $\cN_{\vfiu}$ is dense in $\au$. So we get also that  $\pi(\cN_{\vfiu})$ is dense in $A$ and because $\vfi \neq 0$, we conclude that $\vfiu \neq 0$.
\end{demo}

\medskip

The obvious candidate for the left invariant weight on $\au$ is introduced in the same way. The unitary antipode $\Ru$ introduced in proposition \ref{bi.prop2} will connect both weights.

\begin{proposition}
We define $\vfiu = \vfi \, \pi$. Then $\vfiu$ is a proper weight on $\au$ which has a GNS-construction $(H,\pi,\lau)$ such that $\lau = \la \, \pi$. Moreover, $\vfiu = \psiu \Ru$.
\end{proposition}

Notice that the last equality follows from the commutation $\pi \Ru = R \pi$ and the fact that $\vfi = \psi R$. This last equality also implies that $\vfiu$ is a proper weight.

\bigskip

Since any $^*$-homomorphism sends the open unit ball onto the unit ball of its image (in this case, $A$), the linear mapping  $\pi^* : A^* \rightarrow \au^* : \om \mapsto \om \pi$ is in  an isometry.  Also the Banach space $\au^*$ has a Banach algebra structure with product $\au^* \times \au^* \rightarrow \au^* : (\om,\th) \mapsto \om \th = (\om \ot \th)\deu$. Therefore proposition \ref{com.prop1}.1 implies that $\pi^*$ is an algebra homomorphism.
Also notice that proposition \ref{com.prop3} implies that $\vepu$ is a unit for the Banach algebra $\au^*$.

\begin{proposition} \label{haar.prop1}
The set $\pi^*(A^*)$ is a two-sided ideals in $\au^*$.
\end{proposition}
\begin{demo}
Take $\om \in \au^*$ and $\eta \in A^*$. Then we have for all $x \in \au$ that
$$(\om \, \pi^*(\eta))(x)  =   (\om \ot \pi^*(\eta))(\deu(x))
= (\om \ot \eta)((\io \ot \pi)(\de(x))) = (\om \ot \eta)(\cV^* (1 \ot \pi(x)) \cV) \ , $$
which shows that $\om \, \eta \in \pi^*(A^*)$. So we have proven that  $\pi^*(A^*)$ is a left ideal in $\au^*$. By using the unitary antipode $\Ru$ and the equality $\flip(\Ru \ot \Ru)\deu = \deu \Ru$, we see that $\pi^*(A^*)$ is a two-sided ideal in $\au^*$.
\end{demo}

\medskip

Since $\vfiu = \vfi \pi$, $\psiu = \psi \pi$ and $(\pi \ot \pi)\deu = \de \pi$ propositions 6.2 and 6.3 of \cite{J-V} imply immediately the next result.

\begin{result} \label{haar.res1}
Consider $x \in \au^+$. Then the following holds
\begin{enumerate}
\item If $(\om_{v,v} \ot \io)\deu(x) \in \cM_{\vfiu}^+$ for all $v \in H$, then $x \in \cM_{\vfiu}^+$.
\item If $(\io \ot \om_{v,v})\deu(x) \in \cM_{\psiu}^+$ for all $v \in H$, then $x \in \cM_{\psiu}^+$.
\end{enumerate}
\end{result}

\medskip

Now it is easy to prove the left invariance of $\vfiu$. 

\begin{proposition}
The weight $\vfiu$ is left invariant.
\end{proposition}
\begin{demo}
Choose $x \in \cM_{\vfiu}^+$. Then $\pi(x) \in \cM_{\vfi}^+$. By the left invariance of $\vfi$, we see  for every $\eta \in A^*_+$ that $\pi((\pi^*(\eta) \ot \io)\deu(x)) = (\eta \ot \io)\de(\pi(x)) \in \cM_\vfi^+$ and
$\vfi\bigl(\pi((\pi^*(\eta) \ot \io)\deu(x))\bigr)
= \eta(1) \, \vfi(\pi(x)) = \pi^*(\eta)(1) \,\vfiu(x)$.

So $(\pi^*(\eta) \ot \io)\deu(x)$ belongs to $\cM_{\vfiu}^+$ and 
\begin{equation} \label{haar.eq1}
\vfiu\bigl((\pi^*(\eta) \ot \io)\deu(x)\bigr) = \pi^*(\eta)(1) \,\vfiu(x) \ .
\end{equation}

Now take  $\om \in (\au)^*_+$. Choose $\th \in A^*_+$. By proposition \ref{haar.prop1}, we know that $\om \,\pi^*(\th) \in \pi^*(A^*)$ so the above discussion implies that $(\om \, \pi^*(\th) \ot \io)\deu(x)$ belongs to $\cM_{\vfiu}^+$ and 
$$\vfiu\bigl((\om \,\pi^*(\th)  \ot \io)\deu(x)\bigr) = (\om \,\pi^*(\th))(1) \,\vfiu(x)
= \th(1) \, \om(1) \, \vfiu(x) \ .$$
But $(\pi^*(\th) \ot \io)\deu((\om \ot \io)\deu(x)) = (\om \, \pi^*(\th)  \ot \io)\deu(x)$. 
From result \ref{haar.res1}, we now infer   that $(\om \ot \io)\deu(x) \in \cM_{\vfiu}^+$.

\smallskip

By taking $\th \in A^*_+$ such that $\th(1) = 1$, equation \ref{haar.eq1} gives
$$\om(1) \, \vfiu(x) = \vfiu\bigl((\om \pi^*(\th)  \ot \io)\deu(x)\bigr)
= \vfiu\bigl((\pi^*(\th) \ot \io)\deu((\om \ot \io)\deu(x))\bigr) 
= \vfiu((\om \ot \io)\deu(x))\ .$$
\end{demo}

Since $\psiu = \vfiu \Ru$ and $\flip(\Ru \ot \Ru)\deu = \deu \Ru$, we infer from the previous proposition that 

\begin{corollary}
The weight $\psiu$ is right invariant.
\end{corollary}

\medskip

Although the unitary corepresentation $\cV$ was defined as the generator of a representation of $L^1(\ah)$, it is not so difficult to show that it is the unitary operator naturally associated to the left Haar weight $\vfiu$. We make use of the notations used in section 3.4 of \cite{J-V}.

\begin{proposition}
We have for all $a \in \au$ and $b \in \cN_{\vfiu}$ that 
$\cV^* (a \ot \lau(b)) = (\io \ot \lau)(\deu(b))\,a$.
\end{proposition}
\begin{demo}
Since $\psiu = \psi \pi$, $\psiu$ is easily seen to be approximately KMS. Using proposition 3.21 of \cite{J-V}, we define a unitary  element $V \in \cL(A \ot B_0(H)) = M(A \ot B_0(H))$ such that
$V (a \ot \lau(b)) = (\io \ot \lau)(\deu(b))\,a$ for all $a \in \au$ and $b \in \cN_{\vfiu}$.
We have for all $\om \in \au$, $a,c \in \au$, $b,d \in \cN_{\vfiu}$  that
\begin{eqnarray*}
& & \langle (a \om c^* \ot \io)(V) \lau(b) , \lau(d) \rangle
= \om(\langle V (a \ot \lau(b)) , c \ot \lau(d) \rangle) \\
& & \spat = \om(\langle (\io \ot \lau)(\deu(b))\, a , c \ot \lau(d) \rangle)
= \om\bigl( (\io \ot \vfiu)((c^* \ot d^*)\deu(b)(a \ot 1))\bigr)\\
& & \spat = \vfiu( d^* ( a \om c^* \ot \io)(\deu(b)))
= \langle \lau\bigl(( a \om c^* \ot \io)(\deu(b))\bigr) , \lau(d) \rangle \ ,
\end{eqnarray*}
implying that $(\om \ot \io)(V) \lau(b) = \lau(( \om  \ot \io)(\deu(b)))$ for all $\om \in \au^*$ and every $b \in \cN_{\vfiu}$.  \inlabel{haar.eq2}

Due to the coassociativity of $\deu$, this gives $(\om \ot \io)(V) \, (\th \ot \io)(V) = (\th \om \ot \io)(V)$ for all $\om,\th \in \au^*$. Hence $(\deu \ot \io)(V) = V_{23} V_{13}$.

Equation \ref{haar.eq2} also guarantees that
$(\om \ot \io)((\pi \ot \io)(V)) \, \la(\pi(b)) = \la\bigl((\om \ot \io)\de(\pi(b))\bigr)$ for all $b \in \cN_{\vfiu}$. Since $\pi(\cN_{\vfiu}) = \Nfi$, we infer from result 2.10 of \cite{J-V} that $(\pi \ot \io)(V^*) = W = (\pi \ot \io)(\cV)$. Result \ref{com.res1} allows to conclude that $V = \cV^*$.
\end{demo}

\begin{corollary} \label{haar.cor1}
The slices of $\cV$ are determined by the following formulas:
\begin{enumerate}
\item $(\om \ot \io)(\cV^*) \lau(a) = \lau((\om \ot \io)\deu(a))$ for all $a \in
\cN_{\vfiu}$ and $\om \in \au^*$.
\item $(\io \ot \om_{\la(a),\la(b)})(\cV) = (\io \ot \vfiu)(\deu(b^*)(1 \ot a))$
for all $a,b \in \cN_{\vfiu}$.
\end{enumerate}
\end{corollary}

\bigskip

We will use proposition \ref{bi.prop1} to lift the modular groups of $\vfi$ and $\psi$ to canonical modular groups for $\vfiu$ and $\psiu$ (the canonical nature will be discussed in remark \ref{haar.rem1} and proposition \ref{antipod.prop1}). 

\medskip

We have for $t \in \R$ that $\vfi \sip_t = \nu^t \, \vfi$ and $\vfi \tau_t = \nu^{-t}\,\vfi$.
Define injective positive operators $\nabp$ and $P$ on $H$ such that 
$\nabp^{it} \la(a) = \nu^{-\frac{t}{2}}\, \la(\sip_t(a))$ and 
$P^{it} \la(a) = \nu^{\frac{t}{2}} \,\la(\tau_t(a))$ for all $t \in \R$ and $a \in \Nfi$.
As a matter of fact, $\nabp$ is the modular operator of $\psi$ in the GNS-construction $(H,\io,\ga)$.

\begin{proposition} \label{haar.prop2}
There exists a unique norm continuous one-parameter group $\siup$ on $\au$ such that
$$(\siup_t \ot \io)(\cV) = (1 \ot \nabp^{-it}) \cV (1 \ot P^{-it})$$
for all $t \in \R$. Furthermore, $\pi \siup_t = \si_t \pi$ for all $t \in R$. So $\psiu$ is a KMS weight with modular group $\siup$.
\end{proposition}
\begin{demo}
Let $t \in \R$. By proposition 6.8 of \cite{J-V}, we know that $(\sip_t \ot \tau_{-t})\de = \de \sip_t$.
Therefore the remarks before this proposition and proposition \ref{bi.prop1} imply the existence of a unique $^*$-automorphism $\siup_t$ on $\au$ such that  $(\siup_t \ot \io)(\cV) =(1 \ot \nabp^{-it}) \cV (1 \ot P^{-it})$. By this same proposition, we  know moreover  that
$\pi \siup_t = \si_t \pi$.

So we get for all $\om \in B_0(H)^*$ and $t \in \R$ that 
$$\siup_t((\io \ot \om)(\cV)) = (\io \ot P^{-it} \om \nabp^{-it})(\cV) \ .$$
Therefore equation \ref{gen.eq6} implies easily that the mapping $\R \rightarrow \text{Aut}(\au) : t \mapsto \siup_t$ is a norm continuous one-parameter group on $\au$.
\end{demo}

\begin{definition} \label{haar.def1}
We define the norm continuous one-parameter group $\siu$ on $\au$ such that $\siu_t = \Ru \,\siup_{-t}\, \Ru$ for all $t \in \R$. Then $\pi \siu_t = \si_t \pi$ for all $t \in \R$. So $\vfiu$ is a KMS weight on $\au$ with modular group $\siu$.
\end{definition}

\bigskip

The uniqueness of the Haar weights  on the reduced level implies the uniqueness of the Haar weights on the universal level (see the end of remark 4.4 of \cite{J-V} !)

\begin{theorem} \label{haar.thm1}
Consider a proper weight $\eta$ on $\au$. Then
\begin{itemize}
\item If $\eta$ is left invariant, then there exists a number $r > 0$ such that $\eta = r \, \vfi$.
\item If $\eta$ is right invariant, then there exists a number $r > 0$ such that $\eta = r \, \psi$.
\end{itemize}
\end{theorem}

\bigskip\medskip

Now we formulate the converse of proposition \ref{bi.prop1}. Together with proposition \ref{bi.prop1}, it guarantees the existence of a bijection between bi-automorphisms on $(A,\de)$ and $(\au,\deu)$.

\begin{proposition}
Consider  $^*$-automorphisms $\al$ and $\be$ on $\au$ such that $(\al \ot \be)\deu = \deu \al$. Then the following properties hold.
\begin{enumerate}
\item $(\be \ot \be)\deu = \deu \be$.
\item There exists a number $r > 0$ such that $\vfiu \, \al = r \, \vfiu$ and $\vfiu \, \be = r \, \vfiu$.
\item Define two unitary operators $U$,$V$ on $H$ such that $U \lau(a) = r^{-\frac{1}{2}}\, \lau(\al(a))$ and \newline $V \lau(a) = r^{-\frac{1}{2}}\, \lau(\be(a))$ for all $a \in \cN_{\vfiu}$. Then
$$(\al \ot \io)(\cV) = (1 \ot U^*) \cV (1 \ot V) \hspace{1cm} \text{ and } \hspace{1cm} (\be \ot \io)(\cV)  = (1 \ot V^*) \cV ( 1 \ot V) \ .$$
\item There exists unique $^*$-automorphisms $\alr$,$\ber$ on $A$ such that
$\pi \alr = \pi \al$ and  $\pi \be = \ber \pi$. We have moreover that 
$(\alr \ot \ber)\de= \de \, \alr$.
\end{enumerate}
\end{proposition}

Thanks to theorem \ref{haar.thm1}, result \ref{haar.res1} and corollary \ref{haar.cor1}, the proofs of the first 3 statements are completely analogous as the proofs of the first 3 statements of proposition \ref{bi.prop1}. It is moreover clear that $U \pi(x) U^* = \pi(\al(x))$ and $V \pi(x) V^* = \pi(\be(x))$ for all $x \in \au$, implying that $U A U^* = A$ and $V A V^* = A$. Now define $^*$-automorphisms $\alr$ and $\ber$ such that $\alr(x) = U x U^*$ and $\ber(x) = V x V^*$ for all $x \in A$.

\medskip\medskip

\begin{corollary} \label{haar.cor3}
Consider $^*$-automorphisms $\al_1$,$\al_2$,$\be_1$,$\be_2$ on $\au$ such that $(\al_i \ot \be_i)\deu = \deu \al_i$ (i=1,2). If $\pi \al_1 = \pi \al_2$, then $\al_1 = \al_2$ and $\be_1 = \be_2$.
\end{corollary}
\begin{demo}
Fix $i \in \{1,2\}$. By the previous result, there exists a number $r_i > 0$ such that
$\vfiu \, \al_i =  \vfiu \, \be_i = r_i \, \vfiu$.  Define the unitary operators $U_i$, $V_i$ on $H$ such that $U_i \lau(a) = r_i^{-\frac{1}{2}}\, \lau(\al_i(a))$ and $V_i \lau(a) = r_i^{-\frac{1}{2}}\, \lau(\be_i(a))$ for all $a \in \cN_{\vfiu}$.
We also know that there exist unique $^*$-automorphisms $\al_i',\be_i'$ on $A$ such that
$\al_i' \pi = \pi \al_i$, $\be_i' \pi = \pi \be_i$. It is clear that  $\vfi \, \al_i'  = \vfi \, \be_i' = r_i \, \vfi$ and that $U_i \la(a) = r_i^{-\frac{1}{2}}\, \la(\al_i'(a))$ and $V_i \la(a) = r^{-\frac{1}{2}}\, \la(\be_i'(a))$ for all $a \in \Nfi$.

\smallskip

By assumption $\al_1' = \al_2'$, so $U_1 = U_2$. Also, 
$$(\al_1' \ot \be_1')\de  = \de \al_1' = \de \al_2' = (\al_2' \ot \be_2')\de = (\al_1' \ot \be_2')\de \ ,$$
implying that $(\io \ot \be_1')\de = (\io \ot \be_2')\de$. Therefore the density conditions in definition \ref{red.def1} guarantee that $\be_1'=\be_2'$ and thus $V_1 = V_2$. Consequently, the previous proposition
gives
$$(\al_1 \ot \io)(\cV) = (1 \ot U_1^*) \cV (1 \ot V_1) = (1 \ot U_2^*) \cV (1 \ot V_2)
= (\al_2 \ot \io)(\cV) \ ,$$
thus $\al_1 = \al_2$. Similarly, $\be_1 = \be_2$.
\end{demo}

\medskip

Using the equality $\flip(\Ru \ot \Ru)\deu = \deu \Ru$, we also get the following result.

\begin{corollary} \label{haar.cor2}
Consider $^*$-automorphisms $\al_1$,$\al_2$,$\be_1$,$\be_2$ on $\au$ such that $(\be_i \ot \al_i)\deu = \deu \al_i$ (i=1,2). If $\pi \al_1 = \pi \al_2$, then $\al_1 = \al_2$ and $\be_1 = \be_2$.
\end{corollary}

\medskip

\begin{remark} \rm \label{haar.rem1}
In this remark, we will discuss a first application of this result. Since $\vfiu$ is not faithful in general, the modular group of $\vfiu$ is not uniquely determined. By imposing an extra condition involving the comultiplication, it can be uniquely determined, e.g. in the following way.

\medskip

Consider a norm continuous one-parameter group $\al$ on $\au$ such that
\begin{enumerate} 
\item $\al$ is a modular group for $\vfiu$.
\item For every $t \in \R$, there exists an automorphism $\be_t$ on $\au$ such that $(\be_t \ot \al_t) = \deu \al_t$.
\end{enumerate}
Then $\al$ is equal to $\siu$.

\medskip

Because $\vfi$ is faithful, its modular group $\si$ is uniquely determined. So we get that $\pi \al_t = \si_t \pi = \pi \siu_t$ for all $t \in \R$. But we have also that $(\tauu_t \ot \siu_t)\de = \de \siu_t$ for all $t \in \R$ (see proposition \ref{antipod.prop1}). So corollary \ref{haar.cor2} implies that $\siu = \al$.
\end{remark}

\bigskip\medskip

\sectie{The antipode of the universal quantum group}

In this section, we introduce the polar decomposition through its polar decomposition. The unitary antipode $\Ru$ appeared in proposition \ref{bi.prop2}. It is a $^*$-antiautomorphism on $\au$ such that $\pi \Ru = R \pi$, $\Ru^2 = \io$ and $\flip(\Ru \ot \Ru)\deu = \deu \Ru$.

We first lift the scaling group from the reduced to the universal level and then define  the antipode using both the unitary antipode and the scaling group. In the last part, we easily establish the strong left invariance of $\vfiu$ with respect to our antipode.

\bigskip

Recall that we have for every $t \in \R$ that $(\tau_t \ot \tau_t)\de = \de \tau_t$, $\vfi \tau_t = \nu^{-t} \, \vfi$ and $P^{it} \la(a) = \nu^{\frac{t}{2}}\,\la(\tau_t(a))$ for all $a \in \Nfi$. So proposition \ref{bi.prop1} implies the following one.

\begin{proposition} \label{antipod.prop2}
There exists a unique norm continuous one-parameter group $\tauu$ on $\au$ such that 
$$(\tauu_t \ot \io)(\cV) = (\io \ot P^{-it}) \cV (\io \ot P^{it})$$
for all $t \in \R$.
\end{proposition}

Norm continuity follows in the same way as in  the proof of proposition \ref{haar.prop2}.

\medskip

Referring to proposition 8.23 of \cite{J-V}, we see that $\tauu$ is determined by the fact that $(\tauu_t \ot \hat{\tau}_t)(\cV) = \cV$ for all $t \in \R$. Also recall from proposition \ref{bi.prop2} that $(\Ru \ot \hat{R})(\cV) = \cV$. 

\begin{result} \label{antipod.res1}
We have for all $t \in \R$ that  $\tauu_t \, \Ru = \Ru \, \tauu_t$.
\end{result}
\begin{demo}
By the remarks before proposition 5.22 of \cite{J-V}, we know that $\hat{\tau}_t \hat{R} = \hat{R} \hat{\tau}_t$. Hence
$$(\tauu_t \Ru \ot \hat{\tau}_t \hat{R})(\cV) = \cV = (\Ru \tauu_t \ot \hat{R} \hat{\tau}_t)(\cV) =  (\Ru \tauu_t \ot \hat{\tau}_t \hat{R})(\cV) $$
and the result follows from equation \ref{gen.eq6}
\end{demo}

\begin{proposition} \label{antipod.prop1}
\begin{enumerate}
\item The automorphism groups $\siu,\siup$ and $\tauu$ commute pairwise.
\item We have the following commutation relations for all $t \in \R$:
$$\deu \, \siu_t  =  (\tauu_t \ot \siu_t)\de   \hspace{1.5cm}   \deu\, \siup_t  =  (\siup_t \ot \tauu_{-t}) \deu  \hspace{1.5cm} \deu \,\tauu_t  =  (\tauu_t \ot \tauu_t) \de    $$
\item Let $t \in \R$, then
$$\begin{array}{rclcrcl}
\vfiu \, \siup_t & = & \nu^t \, \vfiu &\hspace{1.5cm}  & \psiu \, \siu_t  &=  &\nu^{-t} \, \psiu  \\
\psiu\, \tauu_t & = & \nu^{-t} \, \psiu & \hspace{1.5cm} & \vfiu \, \tauu_t  &=  &\nu^{-t} \, \vfiu
\end{array}$$
\end{enumerate}
\end{proposition}
\begin{demo}
Let us first comment on the equalities of statement 2. The two last equalities follow from  propositions \ref{bi.prop1}, \ref{haar.prop2} and \ref{antipod.prop2}. The first equality follows from the second equality, definition \ref{haar.def1} and the equality $\flip(\Ru \ot \Ru)\deu$. 
The third statement is an immediate consequence of the corresponding statements on the reduced level (proposition 6.8 of \cite{J-V}).

\smallskip

Since $\tau$ and $\sigma$ commute, the formulas for $\nabp$ and $P$ before proposition \ref{haar.prop2}  imply that the operators $P$ and $\nabp$ commute. So we have for all $s,t \in \R$ that
$$(\tauu_t \,\siup_s \ot \io)\cV = (1 \ot \nabp^{-is} P^{-it}) \cV (1 \ot P^{it} P^{-is})
= (1 \ot  P^{-it} \nabp^{-is}) \cV (1 \ot P^{-is} P^{it})
= (\siup_s \,\tauu_t \ot \io)(\cV) \ ,$$
implying that $\tauu_t \,\siup_s = \siup_s \,\tauu_t$. Combining this with definition \ref{haar.def1} and result \ref{antipod.res1}, we see that $\siu$ and $\tauu$ also commute.

Now the first two equalities in the second statement imply for every $s,t \in \R$ that
$$\deu \siup_s \, \siu_t = (\siup_s \, \tauu_t \ot \tauu_{-s} \, \siu_t) \deu
= (\tauu_t \, \siup_s \ot \siu_t \, \tauu_{-s}) \deu = \deu \siu_t \,\siup_s \ , $$
which by the injectivity of $\deu$ (see the remark after proposition \ref{com.prop3}) gives $\siup_s \,\siu_t = \siu_t \, \siup_s$.
\end{demo}

\medskip

Although both the scaling group $\tauu$ and the unitary antipode $\Ru$ are defined in terms of their behaviour with respect to the unitary corepresentation $\cV$, they can be easily characterized using the projection $\pi : \au \rightarrow A$ and the comultiplication $\deu$.

\begin{proposition}
The following properties characterize $\Ru$ and $\tauu$:
\begin{enumerate}
\item $\tauu$ is the unique norm continuous one-parameter group on $\au$ such that
$\pi \tauu_t = \tau_t \pi$ and $(\tauu_t \ot \tauu_t)\deu = \deu \tauu_t$ for all $t \in \R$.
\item $\Ru$ is the unique $^*$-antiautomorphism on $\au$ such that $\pi \Ru = R \pi$ and $\flip(\Ru \ot \Ru)\deu = \deu \Ru$.
\end{enumerate}
\end{proposition}
\begin{demo}
The statement about $\tauu$ follows immediately from proposition \ref{haar.cor3}. Let us turn to the statement about $\Ru$. So let $\th$ be a $^*$-antiautomorphism on $\au$ such that $\pi \th = R \pi$ and $\flip(\th \ot \th)\deu = \deu \th$. Since $\psi = \vfi R$, we can define an 
anti-unitary operator $U$ on $H$ such that $U \ga(a) =  \la(R(a)^*)$ for all $a \in \Npsi$.
Because $R \pi = \pi \th$, this implies that $\psiu = \vfiu \,\th $ and $U \gau(a) =  \lau(\th(a)^*)$ for all $a \in \cN_{\psiu}$.

\smallskip

Now take $\om \in \au^*$. By corollary \ref{haar.cor1}, we get 
 for all $a \in \cN_{\psiu}$ that 
\begin{eqnarray*}
& & U^* (\om \th^{-1} \ot \io)(\cV)^* U \,\gau(a)
= r^{-\frac{1}{2}}\, U^* (\overline{\om} \th^{-1} \ot \io)(\cV^*) \lau(\th(a)^*)
\\ 
& &  \spat = r^{-\frac{1}{2}}\, U^* \lau\bigl((\overline{\om} \th^{-1} \ot \io)\deu(\th(a)^*)\bigr) = r^{-\frac{1}{2}}\,U^*\lau\bigl(\th((\io \ot \om)(\deu(a))^* )\bigr) \ ,
\end{eqnarray*}
where we used the equality $\flip(\th \ot \th)\deu = \deu \th$ in the last equality. Hence,
$$U^* (\om \th^{-1} \ot \io)(\cV)^* U \, \gau(a) =    \gau((\io \ot \om)\deu(a)) \ .
$$
But this equality has also to hold if $\th = \Ru$, implying that
$$U^* (\om \th^{-1} \ot \io)(\cV))^* U\,\gau(a) = U^* (\om \Ru \ot \io)(\cV))^* U\,\gau(a)$$
for all $a \in \cN_{\psiu}$. 
We conclude from this all that $(\om \th^{-1} \ot \io)(\cV) = (\om \Ru \ot \io)(\cV)$. Therefore equation \ref{gen.eq6} implies that $\om \th^{-1} = \om \Ru$. So we get that $\Ru = \th^{-1}$, thus $\Ru = \th$.
\end{demo}

\bigskip\medskip

As in the case of reduced locally compact quantum groups, the antipode is defined through its polar decomposition.

\begin{proposition}
We define the antipode $\Su = \Ru \tauu_{-\frac{i}{2}} = \tauu_{-\frac{i}{2}} \Ru$. The closed linear mapping $\Su$ satisfies the following basic properties:
\begin{enumerate}
\item $\Su$ is densely defined and has dense range.
\item $\Su$ is injective and $\Su^{-1} = \Ru \, \tauu_{\frac{i}{2}} =
\tauu_{\frac{i}{2}} \, \Ru$.
\item $\Su$ is antimultiplicative : we have for all $x,y \in D(\Su)$ that $xy \in D(\Su)$ and $\Su(xy) = \Su(y) \Su(x)$.
\item We have for all $x \in D(\Su)$ that $\Su(x)^* \in D(\Su)$ and
$\Su(\Su(x)^*)^* = x$.
\item $\Su^2 = \tauu_{-i}$.
\item $\Su \,\Ru = \Ru \,\Su$ and $\Su \,\tauu_t = \tauu_t \, \Su$ for all $t \in \R$.
\end{enumerate}
\end{proposition}

\medskip

Since $\tau_t \pi = \pi \tauu_t$ for all $t \in \R$ and $R \pi = \pi \Ru$, we also get that 
$\pi \Su \subseteq S \pi$. 

\bigskip\medskip

In the next proposition, we show that the weight $\vfiu$ is strongly left invariant with respect to $\Su$. In order to do so, we need to use the Tomita $^*$-algebra $\cT_{\vfiu}$ defined by
$$\cT_{\vfiu} = \{ \, x \in \au \mid x \text{ is analytic with respect to } \siu \text{ and }
\siu_z(x) \in \cN_{\vfiu} \cap \cN_{\vfiu}^* \text{ for all } z \in \C\,\}\ .$$

Let $\nab$ denote the modular conjugation for $\vfiu$ in the GNS-construction $(H,\pi,\lau)$, i.e. $\nab^{it} \lau(a) = \lau(\siu_t(a))$ for all $t \in \R$ and $a \in \cN_{\vfiu}$.
 
For every  $x \in \cT_{\vfiu}$ and $z \in \C$, the element  $\lau(x)$ belongs to $D(\nab^{iz})$ and $\nab^{iz} \lau(x) = \lau(\siu_z(x))$ \ (see e.g. proposition 4.4 of \cite{JK1}).

\begin{proposition} \label{antipod.prop3}
The antipode $\Su$ is characterized by the following properties:
\begin{enumerate}
\item Consider $a,b \in \cN_{\vfiu}$. Then $(\io \ot \vfiu)(\deu(b^*)(1 \ot a))  \in D(\Su)$ and
$$\Su\bigl((\io \ot \vfiu)(\deu(b^*)(1 \ot a))\bigr) = (\io \ot \vfiu)((1 \ot b^*)\de(a))\ .$$
\item The set 
$$\langle \, (\io \ot \vfiu)(\deu(b^*)(1 \ot a)) \mid a,b \in \cN_{\vfiu} \, \rangle$$
is a core for $\Su$.
\end{enumerate}
\end{proposition}
\begin{demo}
\begin{enumerate}
\item Notice that $(\io \ot \vfiu)(\deu(y^*)(1 \ot x)) = (\io \ot \om_{\lau(x),\lau(y)})(\cV)$ for all $x,y \in \cN_{\vfiu}$.

\smallskip

Take $c,d \in \cT_{\vfiu}$. Using proposition \ref{antipod.prop1}, we get for all $t \in \R$ that 
\begin{eqnarray*}
\tauu_t((\io \ot \om_{\lau(c),\lau(d)})(\cV)) & = & \tauu_t\bigl((\io \ot \vfiu)(\deu(d^*)(1 \ot c))\bigr) \\
& = & \tauu_t\bigl((\io \ot \vfiu \siu_t)(\deu(d^*)(1 \ot c))\bigr) \\
& = & (\io \ot \vfiu)\bigl(\deu(\siu_t(d)^*)(1 \ot \siu_t(c))\bigr) \\ 
& = & (\io \ot \om_{\lau(\siu_t(c)),\lau(\siu_t(d))})(\cV) \\ 
& = & (\io \ot \om_{\nab^{it} \lau(c) , \nab^{it} \lau(d)})(\cV) \ .
\end{eqnarray*}
So the remarks before this proposition imply that $(\io \ot \om_{\lau(c),\lau(d)})(\cV)$ belongs to $D(\tauu_{-\frac{i}{2}})$ and 
$$\tauu_{-\frac{i}{2}}((\io \ot \om_{\lau(c),\lau(d)})(\cV))
= (\io \ot \om_{\nab^{\frac{1}{2}} \lau(c) , \nab^{-\frac{1}{2}} \lau(d)})(\cV) \ .$$

Therefore the equality $(\Ru \ot R)(\cV) = \cV$ and the fact that $\hat{R}$ is implemented by $J$ imply that
$(\io \ot \om_{\lau(c),\lau(d)})(\cV)$ belongs to $D(\Su)$ and
\begin{eqnarray*}
& & \Su((\io \ot \om_{\lau(c),\lau(d)})(\cV))
= \Ru((\io \ot \om_{\nab^{\frac{1}{2}} \lau(c) , \nab^{-\frac{1}{2}} \lau(d)})(\cV)) 
= (\io \ot \om_{J \nab^{-\frac{1}{2}} \lau(d) , J \nab^{\frac{1}{2}} \lau(c)})(\cV) \\
& & \spat = (\io \ot \om_{\nab J \nab^{\frac{1}{2}} \lau(d), J \nab^{\frac{1}{2}} \lau(c)})(\cV)
= (\io \ot \om_{\nab \lau(d^*), \lau(c^*)})(\cV)
= (\io \ot \om_{\lau(\siu_{-i}(d^*)) ,\lau(c^*)})(\cV) \\
& & \spat = (\io \ot \vfiu)\bigl(\deu(c)(1 \ot \siu_{-i}(d^*))\bigr)
= (\io \ot \vfiu)((1 \ot d^*)\deu(c))
= (\io \ot \om_{\lau(c),\lau(d)})(\cV^*) \ .
\end{eqnarray*}
Because $\lau(\cT_{\vfiu})$ is dense in $H$, the closedness of $\Su$ implies now easily for every $v,w \in H$ that $(\io \ot \om_{v,w})(\cV)$ belongs to $D(\Su)$ and $\Su((\io \ot \om_{v,w})(\cV)) = (\io \ot \om_{v,w})(\cV^*)$.

If we apply this with $v = \lau(a)$ and $w = \lau(b)$,  statement 1. follows.

\item Proposition \ref{antipod.prop1} and equation \ref{gen.eq6} guarantee that the set under consideration is a dense subset of $D(\tauu_{-\frac{i}{2}})$ which is invariant under $\tauu$. Therefore the result follows e.g. from corollary 1.22 of \cite{JK3}.
\end{enumerate}
\end{demo}

\medskip

Combining this with the equalities $\flip(\Ru \ot \Ru)\deu = \deu \Ru$ and $\vfiu \Ru = \psiu$, we get the following result.

\begin{corollary}
The antipode $\Su$ is characterized by the following properties:
\begin{enumerate}
\item Consider $a,b \in \cN_{\psiu}$. Then $(\psiu \ot \io)((b^* \ot 1)\deu(a))  \in D(\Su)$ and
$$\Su\bigl((\psiu \ot \io)((b^* \ot 1)\deu(a) )\bigr) = (\psiu \ot \io)(\deu(b^*)(a \ot 1))\ .$$
\item The set 
$$\langle \, (\psiu \ot \io)((b^* \ot 1)\deu(a)) \mid a,b \in \cN_{\psiu} \, \rangle$$
is a core for $\Su$.
\end{enumerate}
\end{corollary}

\bigskip

In the proof of proposition \ref{antipod.prop3}.1, we showed that the next  corollary holds for all vector functionals. Since any element in $B_0(H)^*$ can be written as a norm convergent sum of vector functionals, the general result follows immediately from the closedness of $\Su$.

\begin{corollary}
Consider $\om \in B_0(H)^*$. Then $(\io \ot \om)(\cV)$ belongs to $D(\Su)$ and 
$$\Su((\io \ot \om)(\cV)) = (\io \ot \om)(\cV^*) \ .$$
\end{corollary}

\medskip

Also notice that the set $\{\,(\io \ot \om)(\cV) \mid \om \in B_0(H)^*\,\}$ is a core for $\Su$ by proposition \ref{antipod.prop3}.2.

\medskip

\begin{remark} \rm
Let $V$ be a unitary corepresentation of $(\au,\deu)$ on a \cst-algebra $C$. Using the previous corollary (and the remark after it), proposition \ref{com.prop1}.2 and arguing as in the proof of proposition \ref{gen.prop3}, one gets for every $\om \in C^*$ that $(\io \ot \om)(V) \in D(\bar{S}_{\text{\tiny u}})$ and 
$$\bar{S}((\io \ot \om)(V)) = (\io \ot \om)(V^*) \ .$$
Here, $\bar{S}_{\text{\tiny u}}$ denotes the strict closure of $\Su$.
\end{remark}

\bigskip

As in the reduced setting, there exists also a characterization of the antipode solely in terms of the comultiplication. For this, we need some extra terminology (see the beginning of section 5.5 of \cite{J-V}).

\medskip

So fix a \cst-algebra $B$ and an index set $I$. Then we define the following sets:
\begin{enumerate}
\item $MC_I(B) = \{ \, x \text{ an } I\text{-tuple in } M(B) \mid (x_i^*
x_i)_{i \in I} \text{ is strictly summable in } M(B) \, \}$.
\item $MR_I(B) = \{\, x \text{ an } I\text{-tuple in } M(B) \mid (x_i
x_i^*)_{i \in I} \text{ is strictly summable in } M(B) \,\}$.
\end{enumerate}
Both sets are vector spaces under pointwise addition and scalar multiplication. Elements of $MC_I(B)$ can be thought of as infinite columns, elements of $MR_I(B)$ as infinite rows. Notice that the $^*$-operation gives you a bijection between $MC_I(B)$ and $MR_I(B)$.

\medskip

Let $x \in MR_I(B)$ and $y \in MC_I(B)$. Then $(x_i\, y_i)_{i \in I}$ is strictly summable and the net $(\,\sum_{i \in J} x_i \,y_i\,)_{J \in F(I)}$ is bounded.

\medskip

Consider a second \cst-algebra $C$ and  a non-degenerate $^*$-homomorphism $\th$ 
from $B$ into $M(C)$. Then
\begin{enumerate}
\item Let $x \in MR_I(B)$. Then $(\th(x_i))_{i \in I}$ belongs to
$MR_I(C)$.
\item Let $y \in MC_I(B)$. Then $(\th(y_i))_{i \in I}$ belongs to
$MC_I(C)$.
\end{enumerate}

\medskip

Thanks to these comments, the sums appearing in the next two propositions are strictly convergent. The idea of considering elements $a$,$b$ of the form described in the next proposition is due to A. Van Daele.

\begin{proposition} 
Consider $a,b \in \au$ such that there exist an index set $I$, $p \in
MR_I(\au)$ and $q \in MC_I(\au)$ such that
$$a \ot 1 = \sum_{i \in I} \deu(p_i)(1 \ot q_i) \hspace{1.5cm} \text{and} \hspace{1.5cm}
b \ot 1 = \sum_{i \in I} (1 \ot p_i) \deu(q_i).$$
Then $a \in D(\Su)$ and $\Su(a)= b$.
\end{proposition}

\medskip

\begin{proposition}
Define $C$ to be the set consisting of all elements  $a \in \au$ such that there exist an element $b \in \au$, an index set $I$ and $p \in MR_I(\au)$, $q \in MC_I(\au)$ satisfying
$$a \ot 1 = \sum_{i \in I} \deu(p_i)(1 \ot q_i) \hspace{1.5cm} \text{and} \hspace{1.5cm}
b \ot 1 = \sum_{i \in I} (1 \ot p_i) \deu(q_i) \ .$$
Then $C$ is a core for $\Su$.
\end{proposition}

Thanks to the results proven in this section, the proof of the first proposition is completely analogous to the proof of proposition 5.33 of \cite{J-V}.  The proof of the second one is similar to the proof of proposition 5.43 of \cite{J-V}. Using the formula $\flip(\Ru \ot \Ru)\deu = \deu \Ru$, one gets  a variation similar to proposition 5.33 of \cite{J-V}.

\bigskip\medskip

\sectie{The modular element of the universal quantum group}

As explained in section 7 of \cite{J-V}, the modular element $\sde$ of $(A,\de)$ is the unique strictly positive element affiliated with $A$ such that $\si_t(\sde) = \nu^t \, \sde$ for all $t \in \R$ and $\psi = \vfi_\sde$. In this section, we lift $\sde$ to a canonical strictly positive element $\sdeu$ affiliated to $\au$ such that $\sdeu$ is the Radon Nikodym derivative of $\psiu$ with respect to     
$\vfiu$. The basic properties of this modular element of $(\au,\deu)$ are proven.

\medskip

\begin{proposition}
There exists a unique strictly positive element $\sdeu$ affiliated with $\au$ such that
$\pi(\sdeu) = \sde$ and $\deu(\sdeu) = \sdeu \ot \sdeu$. Moreover, $\sdeu \ot \sde = \cV^* (1 \ot \sde)\cV$.
\end{proposition}
\begin{demo}
Let $t \in \R$. We have that $\sde^{it}$ is a unitary element in $M(A)$ such that $\de(\sde^{it}) = \sde^{it} \ot \sde^{it}$. Therefore, proposition \ref{com.prop5} implies  the existence of an element $u_t \in M(\au)$ such that $\pi(u_t) = \sde^{it}$ and $\deu(u_t) = u_t \ot u_t$. Applying $\io \ot \pi$ to this equation and using proposition \ref{com.prop1}.2,
we see that \newline $(1 \ot \sde^{it}) \cV (1 \ot \sde^{-it})= \cV (u_t \ot 1)$. Hence 
$u_t (\io \ot \om)(\cV) = (\io \ot \sde^{-it} \om \sde^{it})(\cV)$
for all $\om \in B_0(H)^*$.

Using equation \ref{gen.eq6}, this implies that the mapping $\R \rightarrow M(\au) :  t \mapsto u_t$ is a strictly continuous  group representation. By the Stone theorem for \cst-algebras (see e.g. theorem 3.10 of \cite{JK3}), we get the existence of a unique strictly positive element $\sdeu$ affiliated  with $\au$ such that $\sdeu^{it} = u_t$ for $t \in \R$. 

Since clearly $\pi(\sdeu)^{it} = \sde^{it}$, $\deu(\sdeu)^{it} = \sdeu^{it} \ot \sdeu^{it}$ and
$\sdeu^{it}  \ot \sde^{it} = \cV^* (1 \ot \sde^{it})\cV$ for all $t \in \R$, the claims of the proposition are satisfied by this $\sdeu$. Uniqueness follows from result \ref{com.res1}.
\end{demo}

The element $\sdeu$ is called the modular element of the quantum group $(\au,\deu)$. We list its basic properties in the next proposition.

\begin{proposition} \label{mod.prop1}
The following properties hold.
\begin{enumerate}
\item $\tauu_t(\sdeu) =  \sdeu$  for $t \in \R$ and $\Ru(\sdeu) = \sdeu^{-1}$.
\item Let $t \in \R$. Then $\sdeu^{it}$ belongs to $D(\bar{S}_{\text{\tiny u}})$ and
$\Su(\sdeu^{it}) = \sdeu^{-it}$.
\item $\siu_t(\sdeu) = \siup_t(\sdeu) = \nu^t \, \sdeu$ for all $t \in \R$.
\item $\siup_t(a) = \sdeu^{it} \, \siu_t(a) \, \sdeu^{-it}$ for all $t\in \R$ and $a \in \au$.
\item $\psiu = (\vfiu)_{\sdeu}$ and $\gau = (\lau)_{\sdeu}$.
\end{enumerate}
\end{proposition}
\begin{demo}
\begin{enumerate}
\item Take $t \in R$. Then $\tauu_t(\sdeu)$ is a strictly positive element such that 
$\pi(\tauu_t(\sdeu)) = \tau_t(\pi(\sdeu)) = \tau_t(\sde) = \sde$ and
$\deu(\tauu_t(\sdeu)) = (\tauu_t \ot \tauu_t)\deu(\sdeu) = \tauu_t(\sdeu) \ot \tauu_t(\sdeu)$. By definition of $\sdeu$, this implies that $\tauu_t(\sdeu) = \sdeu$. Similarly, one proves that
$\Ru(\sdeu)^{-1}  = \sdeu$.
\item Follows immediately from 1.
\item Take $s,t \in \R$. Applying $\siup_t \ot \tauu_t$  to the equation $ \sdeu^{is} \ot \sdeu^{is} = \deu(\sdeu^{is})$ and using proposition \ref{antipod.prop1}.2, we get that
$$   \siup_t(\sdeu^{is}) 
\ot \sdeu^{is} = \siup_t(\sdeu^{is}) \ot \tauu_t(\sdeu^{is}) = \deu(\siup_t(\sde^{is})) \ .$$
If we now apply $\io \ot \pi$ to this equation, proposition \ref{com.prop1}.2 gives
$$\siup_t(\sdeu^{is}) 
\ot \sde^{is} = \cV^* (1 \ot \pi(\siup_t(\sde^{is}))\,) \cV = \cV^*(1 \ot \sip_t(\sde^{is}))\cV
= \nu^{ist} \,\, \cV^*(1 \ot \sde^{is}) \cV = \nu^{ist} \,\, \sdeu^{is} \ot \sde^{is} \ ,$$
implying that $\siup_t(\sdeu^{is}) = \nu^{ist} \, \sdeu^{is}$. 

Combining this with the fact that $\Ru(\sdeu) = \sdeu^{-1}$  and definition \ref{haar.def1}, we also find that $\siu_t(\sdeu^{is}) = \nu^{ist} \, \sdeu^{is}$ for all $s,t \in \R$.
\item Fix $t \in \R$. Define the $^*$-automorphism $\rho$ on $A$ such that $\rho(x) = \sde^{it} \,\tau_t(x)\, \sde^{-it}$ for all $x \in A$. 
By propositions 6.8 and 7.12 of \cite{J-V},  we know that $(\rho \ot \rho)\de = \de \rho$ and $\vfi \rho = \vfi$. Define the unitary operator $U$ on $H$ such that $U \lau(a) = \lau(\rho(a))$ for all $a \in \Nfi$. Then proposition \ref{bi.prop1} implies the existence of a  $^*$-automorphism $\rho_{\text{\tiny u}}$ on $\au$ such that $(\rho_{\text{\tiny u}} \ot \io)(\cV) 
= (1 \ot U^*) \cV (1 \ot U)$. \inlabel{mod.eq1}

Using propositions 6.8 and 7.12 of \cite{J-V}, we get that $(\rho \ot \sip_t)\de = \de \sip_t$ and $(\rho \ot \rho)\de = \de \rho$. Arguing as in the proof of proposition \ref{bi.prop1}, this gives $W(U \ot \nabp^{it}) = (U \ot \nabp^{it}) W$ and $W(U \ot U) = (U \ot U)W$. Hence,  
$(1 \ot U^*) W (1 \ot U) = (U \ot 1) W (U^* \ot 1) = (1 \ot \nabp^{-it}) W (1 \ot \nabp^{it})$. Consequently the definition of $\ah$ implies that $U^* x U = \nabp^{-it} x \nabp^{it}$ for all $x \in \ah$. Since $\cV \in M(\au \ot \ah)$, we conclude from equation \ref{mod.eq1} that $(\rho_{\text{\tiny u}} \ot \io)(\cV) = (1 \ot \nabp^{-it}) \cV (1 \ot \nabp^{it})$.
Remember from proposition \ref{haar.prop2} that $(\siup_t \ot \io)(\cV) = (1 \ot \nabp^{-it}) \cV (1 \ot P^{-it})$. So, arguing as in the proof of proposition \ref{bi.prop1}, we get that
$$((\rho_{\text{\tiny u}} \ot \siup_t)\deu \ot \io)(\cV) = (\deu \siup_t \ot \io)(\cV) \ ,$$
which as usual gives $(\rho_{\text{\tiny u}} \ot \siup_t)\deu = \deu \siup_t$.
Now define $^*$-automorphisms $\al$ and $\be$ on $\au$ such that
$\be(x) = \sde^{-it} \rho_{\text{\tiny u}}(x) \sde^{it}$ and $\al(x) = \sde^{-it} \siup_t(x) \sde^{it}$ for all $x \in \au$. Then the above commutation and the fact that $\deu(\sdeu) = \sdeu \ot \sdeu$ imply that $(\be \ot \al)\deu = \deu \al$. But we have also that $(\tauu_t \ot \siu_t)\deu = \deu \siu_t$ and $\pi \al = \pi \siu_t$ (by proposition 7.12.5 of \cite{J-V}).
Therefore corollary \ref{haar.cor2} guarantees that $\siu_t = \al$.
\item By definition of $\sde$, we  have that $\psi  = \vfi_\sde$. Therefore the next lemma  implies that
$$(\vfiu)_{\sdeu} = \vfi_{\pi(\sdeu)} \, \pi =\vfi_\sde \, \pi = \psi \,
\pi = \psiu \ .$$
The equality $\gau = (\lau)_{\sdeu}$ is proven in a similar way.
\end{enumerate}
\end{demo}

\medskip

\begin{lemma}
Consider a strictly positive element $\al$ affiliated with $\au$ such that there exists a number $r > 0$ such that $\siu_t(\al) = r^t \, \al$ for all $t \in \R$. Then $\si_t(\pi(\al)) = \pi(\al)$ for all $t \in \R$,  $(\vfiu)_\al = \vfi_{\pi(\al)}\,\pi$ and $(\lau)_{\pi(\al)} = \la_{\pi(\al)} \,\pi$.
\end{lemma}
\begin{demo}
We will use the terminology of the first part of section 1.4 of \cite{J-V}. So take an element $x \in \au$ such that $x$ is a left multiplier of $\al^{\frac{1}{2}}$ and $x\,\al^{\frac{1}{2}}$ belongs to $\cN_{\vfiu}$. Because $\vfi \pi = \vfiu$, this implies that $\pi(x\,\al^{\frac{1}{2}}) \in \Nfi$.

Recall that  $\pi(D(\al^{\frac{1}{2}}))\,A$ is a core for $\pi(\al)^{\frac{1}{2}}$ and $\pi(\al)^{\frac{1}{2}}\,(\pi(b)c) = \pi(\al^{\frac{1}{2}}b)\, c$ for all $b \in D(\al^{\frac{1}{2}})$ and $c \in A$. Then it is not so difficult to see that $\pi(x)$ is a left multiplier of $\pi(\al)^{\frac{1}{2}}$ and $\pi(x)\,\pi(\al)^{\frac{1}{2}} = \pi(x \, \al^{\frac{1}{2}}) \in \Nfi$.  

Hence $\pi(x)$ belongs to $\cN_{\vfi_{\pi(\al)}}$ and 
$$\la_{\pi(\al)}(\pi(x)) = \la(\pi(x) \pi(\al)^{\frac{1}{2}}) = \la( \pi(x \, \al^{\frac{1}{2}})) = \lau(x \al^{\frac{1}{2}}) = (\lau)_\al(x) \ .$$
Since such elements $x$ form by definition a core for $(\lau)_\al$, the closedness of $\la_{\pi(\al)}$ implies for every $x \in \cN_{(\vfiu)_\al}$ that $\pi(x) \in \cN_{\vfi_{\pi(\al)}}$ and $\la_{\pi(\al)}(\pi(x)) = (\lau)_\al(x)$. \inlabel{mod.eq2}

As a consequence, $\vfi_{\pi(\al)}\,\pi$ is an extension of $(\vfiu)_\al$. Define the norm continuous one-parameter group $\kappa$ on $\au$ such that $\kappa_t(x) = \al^{it}\, \siu_t(x) \, \al^{-it}$ for all $x \in \au$ and $t \in \R$. Then $\kappa$ is a modular group for $(\vfiu)_\al$.

But we have also that $\pi(\kappa_t(x)) = \pi(\al)^{it} \, \si_t(\pi(x)) \, \pi(\al)^{-it}$ for all $x \in \au$ and $t \in \R$. This implies that $\vfi_{\pi(\al)}\,\pi$ is invariant under $\kappa$. Therefore proposition 1.14 of \cite{JK-Va} guarantees that $\vfi_{\pi(\al)}\,\pi = (\vfiu)_\al$  and thus also $\la_{\pi(\al)} \,\pi = (\lau)_\al$ by equation \ref{mod.eq2}.
\end{demo}

\bigskip

Since the weight $\vfiu$ does not have to be faithful, the Radon Nikodym derivative of $\psiu$ with respect to $\vfiu$ does not have to be unique. By imposing an extra condition, the uniqueness is easily established.

\begin{proposition}
Consider a strictly positive operator $\al$ affiliated with $\au$ such that there exists a number $r > 0$ such that $\siu_t(\al) = r^t \, \al$ for all $t \in \R$ and $\psiu = (\vfiu)_\al$. Then $\al = \sdeu$.
\end{proposition}

Since $(\vfiu)_{\sdeu} = \psiu = (\vfiu)_\al$, we get by proposition 8.41 of \cite{JK1} that $\pi(\al) = \pi(\sdeu) = \sde$. Because  $\deu(\al) = \al \ot \al$, the equality $\al = \sdeu$ holds by the definition of $\sdeu$.

\bigskip

Let us end this section with a last missing commutation relation. Once we have proven this relation, we have shown that $(\au,\deu)$ possesses the same rich analytic structure as $(A,\de)$ does. 

\begin{proposition}
We have for all $t \in R$ that $(\siu_t \ot \siup_{-t})\deu = \deu \tauu_t$.
\end{proposition}
\begin{demo}
Choose $x \in \au$. Propositions \ref{antipod.prop1} and \ref{mod.prop1} imply that 
$$(\siup_t \ot \siu_{-t})(\deu(x)) 
= (\io \ot \siu_{-t} \, \tauu_t)(\deu(\siup_t(x))) 
= (\tauu_t \ot \tauu_t) \bigl(\deu(\siu_{-t}(\siup_t(x)))\bigr)
= \deu(\tauu_t(\sdeu^{it} x \sdeu^{-it})) \ .$$
Using the fact that $\deu(\sdeu) = \sdeu \ot \sdeu$ and proposition \ref{mod.prop1} once more, we infer from the previous chain of equalities that
$$(\siu_t \ot \siup_{-t})(\deu(x))
= (\sdeu^{-it} \ot \sdeu^{-it})\,(\siup_t \ot \siu_{-t})(\deu(x))\,(\sdeu^{it} \ot \sdeu^{it})
= \deu(\tauu_t(x))\ .$$
\end{demo}

\end{document}